\documentclass[12pt]{amsart}
\usepackage{times}
\usepackage{amsfonts}
\usepackage{amsthm}
\usepackage{amsmath}
\advance \topmargin by -\headheight
\advance \topmargin by -\headsep
\evensidemargin \oddsidemargin
\newtheorem{theorem}{Theorem}[section]
\newtheorem{lemma}[theorem]{Lemma}
\newtheorem{corollary}[theorem]{Corollary}

\newtheorem{remark}[theorem]{Remark}

\begin{document}

\title{The free entropy dimension of hyperfinite von Neumann 
algebras}

\author{Kenley Jung}

\dedicatory{For my parents}

\thanks{Research supported in part by the NSF}
\address{Department of Mathematics, University of California,
Berkeley, CA 94720-3840,USA}

\email{factor@math.berkeley.edu}

\begin{abstract} Suppose $M$ is a hyperfinite von Neumann algebra with a
normal, tracial state $\varphi$ and $\{a_1,\ldots,a_n\}$ is a set of
selfadjoint generators for M.  We calculate $\delta_0(a_1,\ldots,a_n)$,
the modified free entropy dimension of $\{a_1,\ldots,a_n\}$.  Moreover we
show that $\delta_0(a_1,\ldots,a_n)$ depends only on $M$ and $\varphi$.  
Consequently $\delta_0(a_1,\ldots,a_n)$ is independent of the choice of
generators for $M$.  In the course of the argument we show that if
$\{b_1,\ldots,b_n\}$ is a set of selfadjoint generators for a von Neumann
algebra $\mathcal R$ with a normal, tracial state and $\{b_1,\ldots,b_n\}$
has finite dimensional approximants, then $\delta_0(N) \leq
\delta_0(b_1,\ldots,b_n)$ for any hyperfinite von Neumann subalgebra $N$
of $\mathcal R.$ Combined with a result by Voiculescu this implies that if
$\mathcal R$ has a regular diffuse hyperfinite von Neumann subalgebra,
then $\delta_0(b_1,\ldots,b_n)=1$.

\end{abstract}

\maketitle

\section{Introduction}

	Suppose $G$ is a group.  Consider the Hilbert space $L^2(G)$ where
$G$ is endowed with counting measure and for each $g \in G$ write $u_g$
for the unitary operator on $L^2(G)$ defined by $(u_g(f))(a) =
f(g^{-1}a)$.  Define the group von Neumann algebra $L(G)$ to be the von
Neumann algebra generated by $\{u_g : g\in G\}$.  It is not too hard to
show that $L(G)$ is a factor (i.e, a von Neumann algebra such that if
$x\in L(G)$ commutes with every other element in $L(G)$, then $x$ is a
scalar multiple of the identity function) iff every nontrivial conjugacy
class of $G$ is infinite.  By definition the free group factor on $m$
generators is $L(F_m)$ where $F_m$ is the free group on $m$ generators.

	Almost two decades ago Dan Voiculescu began to develop a
noncommutative probability theory modeling the free group factors.  The
theory takes the notions of classical probability and transforms them
into ones suited for noncommutative analysis.  Random variables become
elements in von Neumann algebras, expectations turn into normal, tracial
states, and in this particular probability theory, independence always
immediately follows the word `free.' To clarify the last statement
suppose $M$ is a von Neumann algebra with a normal, tracial state
$\varphi$ and $\langle A_j \rangle_{j \in J}$ is a family of unital
$*$-subalgebras of $M$. $\langle A_j \rangle_{j\in J}$ is a freely
independent family provided that for any $j_1,\ldots,j_p \in J$ with
$j_1 \neq j_2,\ldots,j_{p-1} \neq j_p$, and $a_i \in A_{j_i}$

\[ \varphi(a_1) = \cdots = \varphi(a_p)=0 \Rightarrow \varphi(a_1 \cdots 
a_p) = 0. \]     	

\noindent A family of subsets of $M$ is freely independent if the
corresponding family of unital $*$-subalgebras the subsets generate is
freely independent. The definition generalizes the situation in
$L(F_m).$ $L(F_m)$ has a unique normal, faithful, tracial state given by
$\varphi(x)= \langle x\xi, \xi \rangle $ where $\xi$ is the
characteristic function of the identity of $F_m$.  If the generators for
$F_m$ are $g_1,\ldots,g_m,$ then one easily checks that
$\{u_{g_1}\}^{\prime \prime},\ldots,\{u_{g_m}\}^{\prime \prime}$ are
freely independent.

The parallels between classical and free probability go far and the 
interested reader can consult [8] for a general introduction.

	On the operator algebra side, free probability has answered some
open problems in operator algebras.  Developing the ideas of free entropy
and free entropy dimension Voiculescu shows in [10] that the free group
factors possess no Cartan subalgebras (the first known kind with separable
predual).  Ge shows in [3] that the free group factors cannot be
decomposed into a tensor product of two infinite dimensional factors
(again, the first known kind with separable predual) and similarly in [6]
Stefan shows that the free group factors are not the 2-norm closure of the
linear span of a product of abelian $*$-subalgebras.

	In this paper we take a look at what free entropy dimension has in
store for the most tractable kind of von Neumann algebras: those which are
hyperfinite and have a tracial state.  However, free entropy dimension
being the nontrivial machine that it is, we review its definition and
basic properties before stating our results.

\subsection{Definitions and Properties} We recall the concepts of free
entropy and modified free entropy dimension introduced in [9].  For $k, n
\in \mathbb N$ write $M_k^{sa}(\mathbb C)$ for the set of $k \times k$
self-adjoint matrices with complex entries and $(M_k^{sa}(\mathbb C))^n$
for the set of n-tuples of elements in $M_k^{sa}(\mathbb C).$ Suppose
$a_1,\ldots,a_n \in M$ are self-adjoint.  Given $R >0, m,k \in \mathbb N,$
and $\gamma > 0$ define $\Gamma_R(a_1,\ldots,a_n;m,k,\gamma)$ to be the
set of $(x_1,\ldots,x_n)  \in (M_{k}^{sa}(\mathbb C))^n$ such that for
each $j$ $\|x_j\| \leq R$ and for any $1 \leq p \leq m$ and $1 \leq
j_1,\ldots,j_p \leq n$

\[ |tr_k(x_{j_1} \cdots x_{j_p}) - \varphi(a_{j_1} \cdots a_{j_p})| < 
\gamma.\]

\noindent Here $tr_k$ denotes the tracial state on $M_k(\mathbb C)$, the 
$k \times k$ matrices over $\mathbb C$.  If $b_1, \ldots,b_l \in M$, then 
$\Gamma_R(a_1, \ldots, a_n:b_1, \ldots, b_l ; m,k, \gamma)$ denotes the 
set of all $(x_1, \ldots, x_n) \in (M_{k}^{sa}(\mathbb C))^n$ such that 
there exists a $(y_1, \ldots, y_l) \in (M_{k}^{sa}(\mathbb C))^l$ 
satisfying 

\[(x_1, \ldots, x_n, y_1, \ldots, y_l) \in \Gamma_R(a_1, 
\ldots, a_n, b_1, \ldots, b_l ; m, k, \gamma).\]  

For any $d \in \mathbb N$ denote by $\text{vol}$ Lebesgue measure on
$(M_{k}^{sa}(\mathbb C))^d$ (or a subspace thereof) with respect to the
unnormalized Hilbert-Schmidt norm $\|(z_1,\ldots,z_d)\|_2 =
(\sum_{j=1}^d Tr (z_j^2))^{\frac{1}{2}}$ where $Tr$ is the unnormalized
trace.  One successively defines for any $R, \gamma >0$ and $m,k \in
\mathbb N$

\[ \chi_R(a_1,\ldots,a_n:b_1,\ldots,b_l;m,k,\gamma) = k^{-2} \cdot \log 
(\text{vol}(\Gamma_R(a_1,\ldots,a_n:b_1, \ldots, b_l ; m, k, \gamma))) + 
\frac
{n}{2} \log k, \]

\[ \chi_R(a_1, \ldots, a_n: b_1, \ldots, b_l ; m, \gamma) 
=  \limsup_{k \rightarrow \infty} 
\chi_R(a_1,\ldots,z_n:b_1,\ldots,b_l;m,k,\gamma),\]

\[\chi_R(a_1, \ldots, a_n: b_1, \ldots, b_l) = \inf \{ 
\chi_R(a_1,\ldots,a_n:b_1,\ldots,b_l;m,\gamma) : m \in \mathbb N, \gamma 
>0 \}, \]

\[\chi(a_1, \ldots, a_n:b_1, \ldots, b_l) = \sup_{R>0} 
\chi_R(a_1,\ldots,a_n:b_1, \ldots,b_l). \]

\noindent $\chi(a_1,\ldots,a_n:b_1,\ldots,b_l)$ is called the free entropy 
of $a_1,\ldots,a_n$ in the presence of $b_1,\ldots,b_l$.  Replacing the microstate spaces above with $\Gamma_R(a_1,\ldots,a_n;m,k,\gamma)$ yields $\chi(a_1,\ldots,a_n)$ which is simply called the free entropy of $a_1,\ldots,a_n.$

	Now suppose $\{s_1, \ldots, s_n\}$ is a set of freely independent 
semicircular elements in $M$ (by this we mean that 
$\langle \{s_i\}\rangle_{i=1}^{n}$ 
is a family of freely independent sets such that for each 
$1 \leq j \leq n$ $s_j$ is self-adjoint and for any $d \in \mathbb N$ 
$\varphi(s_j^d) = \frac {2}{\pi} \int_{-1}^{1} t^d \sqrt {1-t^2} dt$) such that the von Neumann algebra they generate is freely independent with respect to the strongly closed algebra 
generated by $\{a_1, \ldots, a_n \}$.  Define the modified free entropy dimension of $\{a_1, 
\ldots, a_n \}$ by 

\[ \delta_0(a_1,\ldots,a_n) = n + \limsup_{\epsilon \rightarrow 0} 
\frac {\chi(a_1 + \epsilon s_1, \ldots, a_n + \epsilon s_n : s_1, \ldots, 
s_n)}{| \log \epsilon|}. \]  

	One can view the modified free entropy dimension 
as a noncommutative analogue of the Minkowski dimension.  Given $S 
\subset \mathbb R^d$ the Minkowski dimension of $S$ is $d + \limsup_{\epsilon \rightarrow 0} 
\frac {\mu(N_{\epsilon}(S))}{|\log \epsilon|}$ where $\mu$ is Lebesgue 
measure on $\mathbb R^d$ and $N_{\epsilon}(S)$ is the $\epsilon$-neighborhood
of $S$ (technically this is the upper Minkowski dimension of $S$).  The 
Minkowski dimension defined for $S$ turns out to be the same as the metric 
entropy quantity $\limsup_{\epsilon \rightarrow 0} \frac {\log 
P_{\epsilon}(S)}{|\log \epsilon|}$ where
$P_{\epsilon}(S)$ is the maximum number of points in an $\epsilon$ 
separated subset of $S$ (the $\epsilon$ packing number of $S$).  We will 
have more to say about the connections between Minkowski/metric entropy 
and free entropy dimension.

	Here are a few basic properties of $\delta_0$, all of which are
proven in [10]: \begin{itemize}

\item For $1 \leq j \leq n$ $\delta_0(a_1,\ldots,a_n) \leq
\delta_0(a_1,\ldots,a_j) + \delta_0(a_{j+1},\ldots,a_n)$.  \item For any
$a=a^* \in M$ $\delta_0(a) = 1 - \sum_{t \in sp(a)} (\lambda(\{t\}))^2$
where $\lambda$ is the Borel measure on $sp(a)$ induced by $\varphi$.
\item If $\chi(a_1,\ldots,a_n) > -\infty$, then
$\delta_0(a_1,\ldots,a_n)=n$. \item If $a_1,\ldots,a_n$ are freely
independent, then $\delta_0(a_1,\ldots,a_n) = \delta_0 (a_1) + \ldots
\delta_0(a_n)$.  \end{itemize} 

\noindent Unfortunately, it is not known whether$\delta_0$ is an invariant
of von Neumann algebras with tracial states, i.e., if $\{b_1, \ldots,
b_p\}$ and $\{a_1,\ldots,a_n\}$ are two sets of self-adjoint generators
for $M$(this means that each set generates a strongly closed algebra equal
to $M$ and each element of the set is self-adjoint), then does it follow
that $\delta_0(a_1,\ldots,a_n) = \delta_0 (b_1, \ldots,b_p)$?  An
affirmative answer to this question would show that for $m \neq n$
$L(F_m)$ is not $*$-isomorphic to $L(F_n)$ for it is well known that for
any $m \in \mathbb N$ there exist $m$ semicircular generators
$s_1,\ldots,s_m$ for $L(F_m)$ which satisfy $\delta_0(s_1,\ldots,s_m) =
m$.

\subsection{Results}

	We show that for hyperfinite von Neumann algebras with specified
tracial state $\delta_0$ \emph{is} an invariant.  More specifically,
suppose $M$ is a hyperfinite von Neumann algebra with normal, tracial
state $\varphi.$ By decomposing $M$ over its center it follows that

\begin{eqnarray}& & M \simeq M_0 \oplus \hspace{.05in}
(\oplus_{i=1}^{s} M_{k_i}(\mathbb C)) \oplus M_{\infty} \nonumber \\
& & \varphi \simeq \alpha_0 \varphi_0 \oplus
\hspace{.05in}(\oplus_{i=1}^{s} \alpha_i tr_{k_i}) \oplus 0
\nonumber \end{eqnarray} where $s \in \mathbb N \bigcup \{0\} \bigcup
\{\infty\}, \alpha_i >0$ for $1 \leq i \leq s$ $( i \in \mathbb N),$ $M_0$
is a diffuse von Neumann algebra or $\{0\}$, $\varphi_0$ is a faithful,
tracial state on $M_0$ and $\alpha_0 >0$ if $M_0 \neq \{0\},$ $\varphi_0 =
0$ and $\alpha_0 =0$ if $M_0 = \{0\},$ and $M_{\infty}$ is a von Neumann
algebra or $\{0\}.$ We show that for any self-adjoint generators $a_1,
\ldots,a_n$ for $M$

\[ \delta_0(a_1,\ldots,a_n) = 1 - \sum_{i=1}^{s} \frac
{\alpha_i^2}{k_i^2}. \]

\noindent Because every such $M$ has a finite set of self-adjoint
generators it makes sense to define $\delta_0(M) = \delta_0(a_1, \ldots,
a_n)$ where $\{a_1, \ldots,a_n\}$ is a set of self-adjoint generators for
$M.$ It follows that for any $k \in \mathbb N
\hspace{.1in}\delta_0(M_k(\mathbb C)) = 1 - \frac{1}{k^2}$ and if $M$ is
the hyperfinite $II_1$-factor, then $\delta_0(M) =1$.  The calculations
also show that for hyperfinite $M$ the free entropy dimension number we
obtain for $M$ coincides with the `free dimension' number for $M$ which
appears in Dykema's work [2].

As a consequence of the arguments leading towards the above result we
obtain a `hyperfinite monotonicity' property of $\delta_0$ which says the
following.  Suppose $M$ is an arbitrary von Neumann algebra with specified
tracial state and self-adjoint generators $a_1,\ldots,a_n.$ Assume
moreover that $\{a_1,\ldots,a_n\}$ has finite dimensional approximants,
i.e., for any $m \in \mathbb N,$ $\epsilon >0,$ and $L > \max \{ \|a_i\|
\}_{1 \leq i \leq n}$ there exists an $N \in \mathbb N$ such that for all
$k \geq N$ $\Gamma_L(a_1,\ldots,a_n;m,k,\epsilon) \neq \emptyset.$ If $N$
is a hyperfinite von Neumann subalgebra of $M,$ then

\[ \delta_0(N) \leq \delta_0(a_1,\ldots,a_n).\]
 
\noindent Hyperfinite monotonicity of $\delta_0$ paired with a result by 
Voiculescu ([10]) show that if $M$ has a regular, diffuse, hyperfinite von 
Neumann subalgebra, then $\delta_0(a_1,\ldots,a_n) = 1$ (see Remark 4.8).

	The gist of the argument is simple: Essentially
$\delta_0(a_1,\ldots,a_n)$ is the normalized metric entropy of the unitary
orbit of a well-approximating microstate for $\{a_1,\ldots,a_n\}$.  
Suppose $M$ is hyperfinite with specified tracial state and self-adjoint
generators $\{a_1, \ldots, a_n \}$.  $\chi(a_1 + \epsilon s_1,\ldots,a_n +
\epsilon s_n:s_1, \ldots, s_n)$ is more or less the normalized logarithm
of the volume of the $\epsilon$-neighborhood around the microstates of
$\{a_1, \ldots, a_n \}$.  $M$ being hyperfinite any two such microstates
are approximately unitarily equivalent so $\chi(a_1 + \epsilon s_1,
\ldots, a_n + \epsilon s_n: s_1, \ldots, s_n)$ is a limiting process
calculated from $k^{-2}$ times the logarithm of the volume of the
$\epsilon$-neighborhood of the unitary orbit of a single microstate for
$\{a_1,\ldots,a_n\}$.  Dividing this quantity by $| \log \epsilon |$ and
adding $n$ is close to $k^{-2}$ multiplied by the Minkowski dimension of
the unitary orbit of the microstate or equivalently, the metric entropy of
the set.  Very roughly then, $\delta_0(a_1, \ldots,a_n)$ is the normalized
metric entropy of the unitary orbit of a single well-approximating
microstate for $\{a_1,\ldots,a_n\}.$

	The calculations require more delicacy than we've let on for we
must first fix an $\epsilon$ and find the volume bounds/packing number
bounds with respect to $\epsilon$ not merely over one microstate in one
dimension but over one microstate for each dimension (because the first
process in free entropy takes a limit as the dimensions go to infinity).  
Weak inequalities reduce this to either the investigation of uniform
bounds on the packing numbers of homogeneous spaces obtained from $U_k$,
the $k \times k$ unitaries, or to $\delta_0(a)$ where $a$ is a
self-adjoint element.  In the former case we make crucial use of the
results of Szarek ([7]) and Raymond ([5]).  The latter situation dealing
with $\delta_0(a)$ has already been discussed.

We break up the paper into calculating upper and lower bounds for
$\delta_0(a_1,\ldots,a_n)$ where $a_1,\ldots, a_n$ are arbitrary
self-adjoint generators for hyperfinite $M$ with specified tracial state.  
Section 2 is a short list of notation and assumptions we make throughout
the paper.  Section 3 obtains the upper bound for general $M$.  Section 4
shows that if $\{a_1,\ldots,a_n\}$ is a set of self-adjoint generators of
a diffuse (arbitrary, i.e., not necessarily hyperfinite) von Neumann
algebra with a tracial state and $\{a_1,\ldots,a_n\}$ has finite
dimensional approximants, then $\delta_0(a_1,\ldots,a_n) \geq 1$.  In
particular this yields the desired lower bound for diffuse $M$.  Section 5
obtains the lower bound when $M$ is finite dimensional and Section 6
combines the results of Sections 4 and 5 to arrive at the general lower
bound.  Section 7 gleans immediate corollaries (including hyperfinite
monotonicity of $\delta_0$) and comments on the relation of $\delta_0(M)$
to Dykema's free dimension [2].  Section 8 is an addendum where we prove
some consequences of Szarek's metric entropy bounds of homogeneous spaces.

\section{Definitions and notation}

	Throughout this paper we maintain the notation in the
introduction.  Also we assume throughout that $M$ is a von Neumann algebra
(not necessarily hyperfinite) with separable predual, a unit $I$, and a
normal, tracial state $\varphi$.  $\{s_i:i\in \mathbb N\}$ is always a
semicircular family free with respect to $M$. $\{a_1,\ldots,a_n\} \subset
M$ is a set of self-adjoint generators for $M$ with finite dimensional
approximants.  $R=\max\{\|a_i\|\}_{1\leq i\leq n}.$ Lastly, $|\cdot|_2$ is
the norm on $M_k(\mathbb C)$ or the seminorm on $M$ given by
$|x|_2=(tr_k(x^*x))^{\frac {1}{2}}$ or $|x|_2=(\varphi(x^*x))^{\frac
{1}{2}},$ respectively.

\section{Upper Bound}

	Throughout the section assume that $M$ is hyperfinite and $N$ is 
a finite dimensional $*$-subalgebra of $M$ containing $I.$  Also assume $N$ has self-adjoint generators $\{b_1,\ldots,b_n\}$ 
such that each $b_i$ has operator norm no larger than $R$.  For 
any $k \in \mathbb N$ and $\epsilon>0 \hspace{.1in} T(b_1,\ldots,b_n;k,\epsilon)$ denotes 

\begin{eqnarray} \{(x_1,\ldots,x_n)\in (M_{k}^{sa}(\mathbb C))^n
\hspace{-.2in} & :  \hspace{-.2in} & \text{ there exists a
$*$-homomorphism}\hspace{.1in} \sigma :N \rightarrow M_k(\mathbb
C)\text{such that for all} \nonumber \\ & & 1\leq i\leq n \hspace{.1in}
|\sigma(b_i)-x_i|_{2}\leq \epsilon \text{ and } \|tr_k\circ
\sigma-\varphi|_{N}\|< \epsilon^2\}. \nonumber \end{eqnarray}

	We show that $\delta_0(a_1,\ldots,a_n) \leq 1 - \sum_{i=1}^{s}
\frac {\alpha_i^2}{k_i^2}$ where the $\alpha_i$ and $k_i$ are as in the
canonical decomposition of $M$ discussed on page 3 of the introduction.  
The argument proceeds in several easy steps.  Firstly $\chi(a_1+\epsilon
s_1,\ldots,a_n+\epsilon s_n:s_1,\ldots,s_n)$ is dominated by a number
calculated more or less from $ \text{vol} (T(b_1,\ldots,b_n;k,
\epsilon))$.  Secondly, $T(b_1,\ldots,b_n;k,\epsilon)$ is contained in the
neighborhood of a restricted unitary orbit of any single element in
$T(b_1,\ldots,b_n:k,\epsilon)$.  Szarek's packing number estimates provide
appropriate upper bounds for the volume of the neighborhoods of such
orbits.  Finally by approximating $M$ by fine enough finite-dimensional
$*$-subalgebras $N$ of $M$, standard approximation arguments yield the
promised upper bound.

	The first lemma presented below is standard and we omit the proof.  
It amounts to saying that matricial microstates for self-adjoint 
generators of a finite dimensional von Neumann algebra correspond to 
approximate representations. 

\begin{lemma} For each $\epsilon>0$ there exist an $m\in 
\mathbb N$ and $\gamma>0$ such that for all $k\in \mathbb N$ 

\[ \Gamma_{R+1}(b_1,\ldots,b_n:m,k,\gamma)\subset 
T(b_1,\ldots,b_n;k,\epsilon). \]
\end{lemma}

	Lemma 3.1 has a trivial consequence:

\begin{lemma} Suppose $\epsilon_0>0$ and $\max\{|a_i 
-b_i|_2:1\leq i\leq n\}<\epsilon_0$.  For any $1>\epsilon>0$ 
$\chi(a_1+ \epsilon s_1,\ldots,a_n + \epsilon s_n:s_1,\ldots,s_n)$ is 
dominated by 

\[
\limsup_{k \rightarrow \infty} \left [ k^{-2} \cdot 
\log( \text{vol} (T(b_1,\ldots,b_n;k,3\epsilon+\epsilon_0)))+ \frac{n}{2} 
\cdot \log k \right].
\]

\end{lemma}
\begin{proof}By lemma 3.1 for a given $\epsilon\in (0,1)$ there exist 
$2\leq m \in \mathbb N$ and $\gamma>0$ such that for any $k\in \mathbb N$ 
$\Gamma_{R+1}(b_1,\ldots,b_n;m,k,\gamma)\subset 
T(b_1,\ldots,b_n;k,\epsilon)$.  
For $\gamma^{\prime} >0$ if
\[
(x_1,\ldots,x_n)\in \Gamma_{R+1}(a_1+ \epsilon s_1,\ldots,a_n+ \epsilon 
s_n:b_1,\ldots,b_n,s_1,\ldots,s_n;m,k,\gamma^{\prime}),
\]
then by definition there exists $(y_1,\ldots,y_n,z_1,\ldots,z_n)\in 
(M_{k}^{sa}(\mathbb C))^{2n}$ satisfying 
\[(x_1,\ldots,x_n,y_1,\ldots,y_n,z_1,\ldots,z_n)\in \Gamma_{R+1}(a_1+ 
\epsilon 
s_1,\ldots,a_n+\epsilon 
s_n,b_1,\ldots,b_n,s_1,\ldots,s_n;m,k,\gamma^{\prime}).
\]
By choosing $\gamma^{\prime}< \gamma$ sufficiently small one can force 
$|x_i-y_i|_2< 2\epsilon+ \epsilon_0$ for $1\leq i\leq n$.  Since 
$\gamma^{\prime}< \gamma$ $(y_1,\ldots,y_n)\in 
\Gamma_{R+1}(b_1,\ldots,b_n;m,k,\gamma)$ $\subset 
T(b_1,\ldots,b_n;k,\epsilon)$.  Consequently $(x_1,\ldots,x_n)\in$ 
$T(b_1,\ldots,b_n;k,3\epsilon+ \epsilon_0)$.  We've 
just shown that
\[
\Gamma_{R+1}(a_1+ \epsilon s_1,\ldots,a_n+ \epsilon 
s_n:b_1,\ldots,b_n,s_1,\ldots,s_n;m,k,\gamma^{\prime})\subset 
T(b_1,\ldots,b_n;k,3\epsilon+\epsilon_0).
\]
Basic properties of free entropy imply that $\chi(a_1+ \epsilon s_1,\ldots,a_n+ \epsilon 
s_n:s_1,\ldots,s_n)$ equals

\begin{eqnarray}& & \chi_{R+1}(a_1+ \epsilon 
s_1,\ldots,a_n+ \epsilon s_n:b_1,\ldots,b_n,s_1,\ldots,s_n) \nonumber \\ & 
\leq & \chi_{R+1}(a_1+ \epsilon s_1,\ldots,a_n+ \epsilon 
s_n:b_1,\ldots,b_n,s_1,\ldots,s_n;m,\gamma^{\prime}) \nonumber.
\end{eqnarray}

\noindent By the preceding inclusion the dominating term is less than or
equal to

\[ \limsup_{k \rightarrow \infty} \left [ k^{-2} \cdot \log(
\text{vol}(T(b_1,\ldots,b_n;k,3\epsilon+\epsilon_0)))+ \frac {n}{2} \cdot
\log k \right ]. \] \end{proof}

	Now for the claim which implies that $T(b_1,\ldots,b_n;k, 
\epsilon)$ is contained in the $2\epsilon(1+\sqrt{2}R)$ - neighborhood of  
a restricted unitary orbit of any single element of 
$T(b_1,\ldots,b_n;k,\epsilon)$.

      Because N is finite dimensional and $\varphi$ is tracial assume
from now on that $N \simeq \oplus_{j=1}^p M_{n_j}(\mathbb C)$ and $\varphi|_N \simeq
\oplus_{j=1}^p \alpha_j tr_{n_j}$.

\begin{lemma}For any $k\in\mathbb N$ and $\epsilon>0$ if 
$\sigma_1, \sigma_2:N\rightarrow M_k(\mathbb C)$ are $*$-homomorphisms 
such that $\|tr_k \circ \sigma_1 - tr_k \circ \sigma_2 \|\leq \epsilon^2$, 
then there exists a $u\in U_k$ such that

\[
|u(\sigma_1(x))u^*-\sigma_2(x)|_2\leq 2\|x\|\epsilon.
\] 
\end{lemma}

\begin{proof}Without loss of generality assume that $N=M_{n_1}(\mathbb
C)\oplus \cdots \oplus M_{n_p}(\mathbb C)$.  For any
$l_1,\ldots,l_p\in\mathbb N\bigcup\{0\}$ with $\sum_{i=1}^p n_il_i\leq k$
denote by $\pi_{l_1,\ldots,l_p}:N\rightarrow M_k(\mathbb C)$ the
$*$-homomorphism 

\[ \pi_{l_1,\ldots,l_p}(x_1,\ldots,x_p)= \begin{bmatrix}
x_1\otimes I_{l_1} & 0 & \cdots & 0 \\ 0 & \ddots &  & \vdots \\
\vdots &  & x_p \otimes I_{l_p} & 0 \\ 0 & \cdots & 0 & 0_{p+1} \\
\end{bmatrix} \] 

\noindent where $x_i\otimes I_{l_i}$ is the $n_il_i\times
n_il_i$ matrix with $x_i$ repeated $l_i$ times on the diagonal and
$0_{p+1}$ is the $(k- \sum_{i=1}^p n_il_i) \times (k- \sum_{i=1}^p n_il_i)  
$ 0 matrix.  There exist $m_1,\ldots,m_p,n_1,\ldots,n_p \in \mathbb N \cup
\{0\}$ such that $\sigma_1 \sim \pi_{m_1,\ldots,m_p}, \sigma_2 \sim
\pi_{n_1,\ldots,n_p}$.

	Set $d_i=\min\{m_i,n_i\}$ and observe that if 
$\pi=\pi_{d_1,\ldots,d_p}$, then for $j=1,2\hspace{.1in}\|tr_k\circ \pi - 
tr_k\circ \sigma_j \| \leq \epsilon^2$.  If we can show that for each $j$ there 
exists a $u_j \in U_k$ satisfying 
$|u_j(\sigma_j(x))u_{j}^{*}-\pi (x)|_2\leq 
\|x\|\epsilon$, then we'll be done.  A moment's thought shows that for 
each $j$ there exists a $u_j \in U_k$ such that $Ad u_j \circ \sigma_j - 
\pi$ is a $*$-homomorphism which we'll denote by $\rho_j$.  Obviously 
$\|tr_k \circ \rho_j\| \leq \epsilon^2$ so that for any $x\in N$
\[
|u_j(\sigma_j(x))u_j^*-\pi(x)|_2=((tr_k \circ \rho_j)(x^*x))^{1/2} \leq 
(\epsilon^2 \cdot \|x^*x\|)^{1/2}= \|x\|\epsilon.
\]
\end{proof}

	Given $(x_1,\ldots,x_n), (y_1,\ldots,y_n)\in
T(b_1,\ldots,b_n;k,\epsilon)$ there are representations $\sigma, 
\pi:N \rightarrow M_k(\mathbb C)$ such that for $1 \leq i \leq n
\hspace{.1in}|\sigma (b_i)-x_i|_2, |\pi(b_i)-y_i|_2 < \epsilon$ and
$\|tr_k \circ \sigma - \varphi |_N\|, \|tr_k \circ \pi - \varphi |_N\| <
\epsilon^2$.  So $\|tr_k \circ \sigma - tr_k \circ \pi \| < 2\epsilon^2$.  
By Lemma 3.3 there exists a $u\in U_k$ satisfying \[ |u(\sigma(b_i))u^*-
\pi(b_i)|_2 \leq 2\|b_i\|(\sqrt{2}\epsilon) \leq 2 \sqrt{2} R \epsilon \]
for $1 \leq i \leq n$.  Thus $|ux_iu^*-y_i|_2 \leq 2\epsilon + 2 \sqrt{2}
R \epsilon = 2\epsilon(1+\sqrt{2}R)$.  From now on for $z\in
M_{k}^{sa}(\mathbb C)$ and $\gamma>0$ define $B(z,\gamma)=\{x \in
M_{k}^{sa}(\mathbb C): |x-z|_2<\gamma\}$.  We've just proved:

\begin{corollary}If $(x_1,\ldots,x_n) \in T(b_1,\ldots,b_n;k,\epsilon)$, 
then $T(b_1,\ldots,b_n;k,\epsilon)$ is contained in 
\[ \bigcup_{u \in U_k} 
[B(ux_1u^*,2\epsilon(1+\sqrt{2}R)) \times \cdots \times 
B(ux_nu^*,2\epsilon(1+\sqrt{2}R))].\] 
\end{corollary}

	We remark here that Lemma 3.1, Lemma 3.3, and Corollary 3.4 also holds in the situation where $N$ does not contain $I$.    

	We now draw out a trivial consequence of Szarek's estimates for 
covering numbers of homogeneous spaces.  These results are the heart of 
the calculation of the upper bound.

	For $k\in \mathbb N$ suppose $m,k_1,\ldots,k_m,l_1,\ldots,l_m \in 
\mathbb N, \sum_{i=1}^{m} k_il_i=k$, and $H \subset U_k$ is a proper 
Lie subgroup of $U_k$ consisting of all matrices of the form
\[
\begin{bmatrix}
u_1\otimes I_{k_1} & \cdots & 0 \\
\vdots  & \ddots & \vdots \\
0  & \cdots & u_m \otimes I_{k_m} \\
\end{bmatrix}
\]
where $u_i\in U_{l_i}$ and $u_i \otimes I_{k_i}$ is the $k_il_i \times 
k_il_i$ matrix obtained by repeating $u_i$ $k_i$ times along the diagonal.  
Such Lie subgroups $H$ of $U_k$ will be called tractable.

	A simple application of Theorem 11 in [7] yields:
\begin{lemma} There exist constants $C, \beta > 0$ such 
that for any $k\in \mathbb N$, any tractable Lie subgroup $H$ of $U_k$ and 
$\epsilon \in (0, \beta )$,
\[
N(X,\epsilon) \leq \left(\frac{C}{\epsilon}\right)^d
\]
where X is the manifold $U_k/H$ endowed with the quotient metric induced 
by the operator norm, $d$ is the real dimension of $X$, and 
$N(X,\epsilon)$ is the minimum number of balls of radius $\epsilon$ 
required to cover $X$.
\end{lemma}

	We sequester a rigorous demonstration to the Addendum.

\begin{lemma}If $1>r>0,$ then there exists a $k_0 \in \mathbb N$ such that
for each $k>k_0$ there is a corresponding $*$-homomorphism $\sigma_k:N
\rightarrow M_k(\mathbb C)$ satisfying: \begin{itemize} \item $\|tr_k
\circ \sigma_k - \varphi|_N \|< r^2.$ \item The unitaries $H_k$ of
$\sigma_k(N)^{\prime}$ is a tractable Lie subgroup of $U_k$ satisfying \[
k^2 \left (1- r - \sum_{i=1}^p \frac {\alpha_{i}^{2}}{n_{i}^{2}} \right)
<\dim(U_k/H_k)< k^2 \left (1 + r - \sum_{i=1}^{p} \frac
{\alpha_{i}^2}{n_i^2} \right ).  \] \end{itemize} \end{lemma}

\begin{proof}Suppose $1 > \varepsilon >0.$ Choose $n_0 \in \mathbb N$ such
that $\frac {1}{n_0} < \frac {\varepsilon}{p+1}.$ Set
$k_0=(n_0+1)n_1\cdots n_p$.  Suppose $k>k_0.$ Find the unique $n \in
\mathbb N$ (dependent on $k$)  satisfying \[ nn_1 \cdots n_p \leq k <
(n+1)n_1 \cdots n_p. \] Set $d=n n_1 \cdots n_p$ ($d$ dependent on $k$)and
find $m_1,\ldots,m_p \in \mathbb N \bigcup \{0\}$ satisfying $\alpha_i -
\varepsilon < \frac {m_i}{n} < \alpha_i + \varepsilon$ for each $i$ and
$\sum_{i=1}^{p} \frac {m_i}{n} =1$ (the $m_i$ depend on $k$).  Set
$l_i(k)= \frac {d m_i}{n n_i} \in \mathbb N \bigcup \{0\}$ and $l_{p+1}(k)
= k - \sum_{i=1}^p l_i(k)n_i.$ Assume without loss of generality that
$N=M_{n_1}(\mathbb C)\oplus \cdots \oplus M_{n_p}(\mathbb C)$ and
$\varphi|_N = \alpha_1 tr_{n_1} \oplus \cdots\oplus \alpha_p
tr_{n_p}$.  Define $\sigma_k:N \rightarrow M_k(\mathbb C)$ by 

\[\sigma_k(x_1,\ldots,x_n)= \begin{bmatrix} I_{l_1(k)} \otimes x_1 & 0 &
\cdots & 0 \\ 0 & \ddots &  & \vdots \\ \vdots &  & I_{l_p(k)}
\otimes x_p & 0 \\ 0 & \cdots & 0 & 0_{l_{p+1}(k)} \\ 
\end{bmatrix} \]

\noindent where $0_{l_{p+1}(k)}$ is the $l_{p+1}(k) \times l_{p+1}(k)$ 0
matrix and $I_{l_i(k)} \otimes x_i$ is the $l_i(k) n_i \times
l_i(k)n_i$ matrix obtained by taking each entry of $x_i, (x_i)_{st}$, and
stretching it out into $(x_i)_{st} \cdot I_{l_i(k)}$ where $I_{l_i(k)}$ is
the $l_i(k) \times l_i(k)$ identity matrix. \[ (tr_k \circ
\sigma)(x_1,\ldots,x_p) = \frac {1}{k} \cdot \sum_{i=1}^p l_i(k)  \cdot
Tr(x_i) = \sum_{i=1}^{p} \frac {d m_i}{k n } \cdot tr_{n_i}(x_i). \]
$\frac {d}{k} > 1- \varepsilon$ so $\alpha_i + \varepsilon \geq \frac
{d}{k} \cdot \frac {m_i}{n} > (\alpha_i - \varepsilon)(1- \varepsilon) >
\alpha_i - 2 \varepsilon$.  It follows that $\|tr_k \circ \sigma_k -
\varphi|_N\| < 2 p \varepsilon $.

	$H_k$ consists of all matrices of the form \[ \begin{bmatrix}
u_1\otimes I_{n_1} & 0 & \cdots & 0 \\ 0 & \ddots & & \vdots \\ \vdots &
& u_p \otimes I_{n_p} & 0 \\ 0 & \cdots & 0 & u_{p+1} \\ \end{bmatrix}
\] where $u_i \in U_{l_i(k)}$ for $1 \leq i \leq p+1$ and $u_i \otimes
I_{n_i}$ is the $l_i(k) n_i \times l_i(k) n_i$ matrix obtained by
repeating $u_i \hspace{.1in} n_i$ times along the diagonal.  $H_k$ is
obviously tractable.

  For a lower bound on $\dim(U_k/H_k)$ we have the estimate:
\[
l_{p+1}(k)= k - \sum_{i=1}^{p} \frac{d m_i}{n} = k - d < n_1 \cdots n_p < \frac{k_0}{n_0} < k \cdot \varepsilon
\]
so that
\[ \dim H_k = l_{p+1}(k)^2 + \sum_{i=1}^{p}l_i(k)^{2} < k^2 \left 
(\varepsilon 
+ \sum_{i=1}^{p} \frac {(\alpha_i + \varepsilon)^2}{n_i^2} \right) < k^2\left ( 3p\varepsilon + \sum_{i=1}^{p}\frac {\alpha_{i}^{2}}{n_{i}^2} 
\right ). 
\]
\noindent Hence, $\dim (U_k/H_k) = k^2 - \dim H_k$ is bounded from below by $k^2(1 - 3p\varepsilon - \sum_{i=1}^p \frac {\alpha_i^2}{n_i^2}).$

For an upper 
bound on $\dim (U_k/H_k)$ observe that $\dim H_k>\sum_{i=1}^{p}l_i(k)^{2}$ 
whence
\begin{eqnarray}\dim(U_k/H_k) = k^2 - \dim H_k < k^2 - 
\sum_{i=1}^{p}l_i(k)^{2} & 
< & k^2 \left (1-\sum_{i=1}^{p} \frac {(\alpha_i - 2\varepsilon)^2}{n_i^2} \right) 
\nonumber \\ & < & k^2 \left ( 1 + 4 \varepsilon - \sum_{i=1}^p \frac {\alpha_i^2}{n_i^2} \right). \nonumber
\end{eqnarray}

Set $\varepsilon = \frac {r^2}{4p}.$
\end{proof}

	We now make the key calculation on the upper bound of lemma 3.2.

\begin{lemma}For $\min\{\beta,C\}>\epsilon>0$ 
\[\limsup_{k\rightarrow \infty}\hspace{.1in}\left [k^{-2}\cdot \log( 
\text{vol} (T(b_1,\ldots,b_n;k, \epsilon))) + \frac {n}{2} \cdot \log k 
\right] \leq \log(\epsilon^{n-\bigtriangleup}) + \log D
\]
where $\bigtriangleup= 1- \sum_{i=1}^{p} \frac 
{\alpha_{i}^{2}}{n_{i}^{2}}$ and $D
= \pi^{\frac {n}{2}}(8(R+1))^nC^{\bigtriangleup}[(2e)^{\frac {n}{2}}]$.
\end{lemma}

\begin{proof}Suppose $\min\{\beta,C\}>\epsilon>r>0$.  By lemma 3.6 
there is a $k_0\in \mathbb N$ such that for each $k \geq k_0$ there exists a 
$*$-homomorphism $\sigma_k: N \rightarrow M_k(\mathbb C)$ satisfying $\|tr_k 
\circ \sigma_k - \varphi|_N \| < r^2$ and the additional condition that if 
$H_k$ is the unitary group of $\sigma_k(N)^{\prime}$, then $H_k$ is 
tractable and
\[\dim(U_k/H_k)< 
k^2 \left (1 + r - \sum_{i=1}^{p} \frac {\alpha_{i}^2}{n_i^2} \right ).\] 
Set $d_k=\dim(U_k/H_k)$ and 
$m_r= -r + \sum_{i=1}^{p} \frac {\alpha_{i}^2}{n_i^2} $.  There exists a 
set $<u_{k,s}>_{s\in S_k}$ contained in $U_k$ such that for each $u\in 
U_k$ there exists an $s\in S_k$ and $h\in H_k$ satisfying 
$\|u-u_{k,s}h\|<\epsilon$ and $S_k$ has cardinality not exceeding 
$\left ( \frac {C}{\epsilon} \right)^{d_k} < \left (\frac {C}{\epsilon} 
\right )^{k^2(1 - m_r)}$.

	Set $L=8(R+1)$.  I claim that
\[T(b_1,\ldots,b_n;k,\epsilon) \subset \bigcup_{s\in S_k} 
[B(u_{k,s}\sigma_k(b_1)u_{k,s}^{*},L \epsilon) \times \cdots 
\times B(u_{k,s}\sigma_k(b_n)u_{k,s}^{*},L \epsilon)]. \]
Suppose $(x_1,\ldots,x_n)\in T(b_1,\ldots,b_n;k,\epsilon).$ Clearly $(\sigma_k(b_1),\ldots,\sigma_k(b_n)) \in 
T(b_1,\ldots,b_n;k,\epsilon).$  By Corollary 3.4 there exists a 
$u\in U_k$ such that for all $1 \leq i \leq n \hspace{.1in} |x_i - 
u\sigma_k(b_i)u^{*}|_2 \leq 2\epsilon(1+\sqrt{2}R)$.  There 
exists an $s\in S_k$ and $h\in H_k$ such that $\|u-u_{k,s}h\| < \epsilon$.  
Now
\begin{eqnarray}
|u\sigma_k(b_i)u^*-u_{k,s}\sigma_k(b_i)u_{k,s}^{*}|_2 & = & 
|u\sigma_k(b_i)u^*-u_{k,s}h\sigma_k(b_i)h^*u_{k,s}^{*}|_2 \nonumber \\ & 
\leq & \|u-u_{k,s}h\| \cdot |\sigma_k(b_i)u^*|_2 + 
|u_{k,s}h\sigma_k(b_i)|_2 \cdot \|u^*-h^*u_{k,s}^{*}\| \nonumber \\ & \leq 
& 2\epsilon R. \nonumber\end{eqnarray}
Hence for all $1 \leq i \leq n \hspace{.1in} 
|x_i-u_{k,s}\sigma_k(b_i)u_{k,s}^{*}|_2 
\leq 2\epsilon (1 + \sqrt{2}R) + 2 \epsilon R< L \epsilon$. 

	By the inclusion demonstrated in the preceding paragraph $\log 
(\text{vol}(T(b_1,\ldots,b_n;k,\epsilon)))$ is dominated by
\begin{eqnarray} \log \left (|S_k| \cdot \frac {\pi^{\frac 
{nk^2}{2}}(L \sqrt{k} 
\epsilon)^{nk^2}}{(\Gamma(\frac{k^2}{2} +1))^n} \right) & = & 
\log \left (|S_k| \cdot \frac {(\pi^{\frac 
{n}{2}} L^n k^{\frac {n}{2}} \epsilon^n)^{k^2}}{(\Gamma(\frac{k^2}{2} 
+1))^n} \right) \nonumber \\ & \leq 
& k^2 \cdot \log \left[ \left (\frac {C}{\epsilon} 
\right)^{1-m_r}(\pi^{\frac {n}{2}} L^n 
k^{\frac {n}{2}} \epsilon^n) \right ] - n \cdot \log \Gamma \left (\frac 
{k^2}{2} 
+1 \right) \nonumber \\ & \leq & k^2 \cdot \log(\pi^{\frac {n}{2}} L^n 
C^{1-m_r} k^{\frac {n}{2}} \epsilon^{n-1+m_r}) - n \cdot \log \left[ 
\left( \frac {k^2}{2e} \right )^{\frac {k^2}{2}} \right] \nonumber \\ & = 
& k^2 \cdot \log(\pi^{\frac {n}{2}} L^n
C^{1-m_r} \epsilon^{n-1+m_r}) - \frac{nk^2}{2} \cdot \log k 
+k^2\log[(2e)^{\frac{n}{2}}] \nonumber\end{eqnarray} 
provided $k \geq k_0$.  $\limsup_{k \rightarrow \infty} 
\left [ k^{-2}\log(\text{vol}(T(b_1,\ldots,b_n;k, \epsilon))) + 
\frac{n}{2} \cdot \log k \right]$ 
is therefore dominated by 
\[ \limsup_{k \rightarrow \infty} (\log(\pi^{\frac{n}{2}} L^n C^{1-m_r} 
\epsilon^{n-1+m_r}) + \log[(2e)^{\frac{n}{2}}] = 
\log(\epsilon^{n-1+m_r}) + \log (\pi^{\frac {n}{2}} L^n 
C^{1-m_r}[(2e)^{\frac {n}{2}}]).\]
Hence
\begin{eqnarray}\log(\epsilon^{n-\bigtriangleup}) + \log D & = & \log 
(\epsilon^{n-\bigtriangleup}) + \log(\pi^{\frac {n}{2}} L^n 
C^{\bigtriangleup} 
[(2e)^{\frac{n}{2}}]) \nonumber \\ & = & \lim_{r \rightarrow 
0}\hspace{.1in} 
[(\log(\epsilon^{n-1+m_r}) + \log(\pi^{\frac {n}{2}} L^n 
C^{1-m_r}[(2e)^{\frac{n}{2}}])] \nonumber \\ & \geq & \limsup_{k 
\rightarrow \infty}\hspace{.1in} 
\left [k^{-2} \cdot \log(\text{vol}(T(b_1,\ldots,b_n;k,\epsilon))) + 
\frac{n}{2} 
\cdot 
\log k \right]. 
\nonumber \end{eqnarray}
\end{proof}

	If $M$ is finite dimensional, then lemma 3.7 yields the desired 
upper bound for $\delta_0(a_1,\ldots,a_n)$.  With just a few more easy 
observations Lemma 3.7 allows us to bootstrap the upper bound for 
$\delta_0(a_1,\ldots,a_m)$ in the general situation.

	If $B$ is a finite dimensional von Neumann algebra with a 
positive trace $\psi$ and $B \simeq \oplus_{i=1}^{s} M_{q_i}(\mathbb 
C),$ $\psi \simeq \oplus_{i=1}^{s} r_i tr_{q_i}$, define 
$\bigtriangleup_{\psi}(B) = 1-\sum_{i=1}^{s} \frac {r_{i}^{2}}{q_{i}^{2}}$.  
Clearly $\bigtriangleup_{\psi}(B)$ is well-defined.

\begin{lemma}If $A\subset B$ is a unital inclusion of finite dimensional
von Neumann algebras, and $\psi$ is a positive trace on $B$, then
$\bigtriangleup_{\psi}(A) \leq \bigtriangleup_{\psi}(B)$. \end{lemma}

\begin{proof}By assumption $B \simeq \oplus_{j=1}^{s} M_{q_j}(\mathbb
C)$ and $\psi \simeq \oplus_{j=1}^{s} r_j tr_{q_j}$ for some
$s,q_1,\ldots,q_s,r_1,\ldots,r_s \in \mathbb N$.  A is $*$-isomorphic to
$\oplus_{i=1}^{d} M_{p_i}(\mathbb C)$ for some $d, p_1,\ldots,p_d \in
\mathbb N$.  Denote $\langle \Lambda_{ij} \rangle_{1\leq i \leq d, 1\leq
j \leq s}$ to be the inclusion matrix of $A$ into $B$ with respect to
the dimension vectors $\langle p_i \rangle_{i=1}^{d}$ and $\langle q_j
\rangle_{j=1}^{s}$ for $A$ and $B$, respectively.  Since $A \subset B$
is a unital inclusion

\[\sum_{i=1}^{d} \sum_{j=1}^{s} \frac 
{(\Lambda_{ij}p_ir_j)^2}{q_j^2} \cdot \frac {1}{p_{i}^{2}} = 
\sum_{j=1}^{s} 
\sum_{i=1}^{d} \frac {\Lambda_{ij}^{2} r_{j}^{2}}{q_{j}^{2}} \geq  
\sum_{j=1}^{s} \frac {r_{j}^{2}}{q_{j}^{2}}. \]

\noindent $\bigtriangleup_{\psi}(A) = 1 - \sum_{i=1}^{d} \left (\left 
(\sum_{j=1}^{s} \frac {\Lambda_{ij}p_ir_j}{q_j} \right )^2 \cdot \frac 
{1}{p_{i}^{2}} 
\right ) 
\leq 
1 - \sum_{j=1}^{s} 
\frac {r_{j}^{2}}{q_{j}^2} = \bigtriangleup_{\psi}(B).$\end{proof} 

\begin{lemma}$\delta_0(a_1,\ldots,a_n)=\delta_0(a_1,\ldots,a_n,I).$\end{lemma}
\begin{proof}By Propositions 6.4 and 6.6 of [10] and Proposition 6.3 of 
[9]
\[\delta_0(a_1,\ldots,a_n,I) \leq \delta_0(a_1,\ldots,a_n) + \delta_0(I) = \delta_0(a_1,\ldots,a_n).\]

On the other hand since the strongly closed $*$-algebra generated by 
$\{a_1,\ldots,a_n\}$ is $M$, by Theorem 4.3 of [10] 
$\delta_0(a_1,\ldots,a_n)\leq \delta_0(a_1,\ldots,a_n,I)$.\end{proof}

	We're now in a position to calculate the upper bound for
$\delta_0(a_1,\ldots,a_n)$.  By decomposing $M$ over its center it
follows that \[M \simeq M_0 \oplus \hspace{.05in} (\oplus_{i=1}^{s}
M_{k_i}(\mathbb C)) \oplus M_{\infty}, \hspace{.15in} \varphi \simeq
\alpha_0 \varphi_0 \oplus \hspace{.05in}(\oplus_{i=1}^{s} \alpha_i
tr_{k_i}) \oplus 0 \] where all quantities above are as in the
introduction. Write $I_i$ for the identity of $M_{k_i}(\mathbb C)$ for
$1 \leq i \leq s$ and if $M_0 \neq \{0\}$, then write $I_0$ for the
identity of $M_0$.  A moment's thought shows that for the purposes of
the theorem below we can neglect the $M_{\infty}$ summand and assume \[M
= M_0 \oplus \hspace{.05in} (\oplus_{i=1}^{s} M_{k_i}(\mathbb C)),
\hspace{.1in} \varphi = \alpha_0 \varphi_0 \oplus
\hspace{.05in}(\oplus_{i=1}^{s} \alpha_i tr_{k_i}).\]
\begin{theorem}$\delta_0(a_1,\ldots,a_n) \leq 1 - \sum_{i=1}^{s} \frac
{\alpha_{i}^{2}}{k_{i}^{2}}$. \end{theorem} \begin{proof}Set $\alpha=1 -
\sum_{i=1}^{s} \frac{\alpha_{i}^{2}}{k_{i}^{2}}$.  There exists a nested
sequence of finite dimensional $*$-subalgebras of $M$, $\langle N_m
\rangle _{m=1}^{\infty}$, such that $\bigcup_{m=1}^{\infty}N_m$ is
strongly dense in $M$, for each $m\in \mathbb N \hspace{.1in} I\in N_m$,
and $\lim_{m \rightarrow \infty} \bigtriangleup_{\varphi}(N_m) \leq
\alpha$.  This is clear for if $M_0=\{0\}$, then for each $m$ define
\[N_m=0 \oplus (\oplus_{1 \leq j \leq \min\{m,s\} } M_{k_j}(\mathbb C))
\oplus \mathbb C \cdot (\oplus_{m < j \leq s} I_j).\] Observe that
\[\bigtriangleup_{\varphi}(N_m)=1- \sum_{j=1}^{\min \{m,s\}} \frac
{\alpha_{j}^{2}}{k_{j}^{2}} - (\sum_{m< j \leq s} \alpha_j)^2.\]
$\lim_{m\rightarrow \infty} \bigtriangleup_{\varphi}(N_m) = \alpha$ and
all the other properties required of the $N_m$ are easily checked.  If
$M_0 \neq \{0\}$, then there exists a nested sequence of finite
dimensional $*$-subalgebras of $M_0$, $\langle A_m
\rangle_{m=1}^{\infty}$ with $I_0\in A_m$ for each $m$ and
$\bigcup_{m=1}^{\infty}A_m$ strongly dense in $M_0$.  For each $m$
define \[N_m=A_m \oplus (\oplus_{1 \leq j \leq
\min\{m,s\}} M_{k_j}(\mathbb C)) \oplus \mathbb C \cdot
(\oplus_{m < j \leq s} I_j).\] Observe that
\begin{eqnarray}\bigtriangleup_{\varphi}(N_m) & = & 1 +
(\bigtriangleup_{\alpha_0\varphi_0}(A_m) - 1) \hspace{.05in} -
\hspace{-.1in} \sum_{j=1}^{\min \{m,s\}} \alpha_{j}^{2}/k_{j}^{2}
\hspace{.05in} - \hspace{.05in} (\sum_{m < j \leq s } \alpha_j)^2
\nonumber \\ & \leq & 1 - \hspace{-.1in} \sum_{j=1}^{\min \{m,s\}}
\hspace{-.1in} \alpha_{j}^{2}/k_{j}^{2} \hspace{.05in} - \hspace{.05in}
(\sum_{m < j \leq s} \alpha_j)^2. \nonumber \end{eqnarray}

As $m \rightarrow \infty$ the dominating term above converges to $\alpha$ 
so $\lim_{m\rightarrow \infty} \bigtriangleup_{\varphi}(N_m) \leq \alpha$ 
(existence of the limit is guaranteed by Lemma 3.8 and the fact that $N_m 
\subset N_{m+1}$).  All the other properties of the $N_m$ are easily 
checked.  Notice that in either cases $\lim_{m \rightarrow \infty} 
\bigtriangleup_{\varphi}(N_m) \leq \alpha$ and Lemma 
3.8 imply $\bigtriangleup_{\varphi}(N_m) \leq \alpha$ for all $m\in 
\mathbb N$.

	Take a sequence $\langle N_m \rangle_{m=1}^{\infty}$ as
constructed above.  Suppose $\min\{1,\beta/4,C\} > \epsilon > 0$.  By
Kaplansky's Density Theorem there exists an $m_0 \in \mathbb N$ and
self-adjoint $x_1,\ldots,x_n \in N_{m_0}$ satisfying $|x_i -a_i|_2 <
\epsilon$ and $\|x_i \| \leq R$ for $1 \leq i \leq n$.  Denote by $B$
the $*$-algebra generated by $\{x_1,\ldots,x_n, I \}$.  By Lemma 3.2 and
Lemma 3.7 $\chi(a_1+ \epsilon s_1,\ldots,a_n + \epsilon s_n, I+ \epsilon
s_{n+1}:  s_1,\ldots,s_{n+1})$ is dominated by

\begin{eqnarray*}\limsup_{k \rightarrow \infty} \left [k^{-2} \cdot
\log(\text{vol}(T(x_1,\ldots,x_n,I;k, 4\epsilon))) + \frac {n+1}{2} \cdot
\log k \right ] & \leq &
\log((4\epsilon)^{n+1-\bigtriangleup_{\varphi}(B)}) + \log D \\ & \leq &
\log (\epsilon^{n+1- \bigtriangleup_{\varphi}(B)}) + \\ & & \log
(4^{n+1}D) \\ \end{eqnarray*}

\noindent where $D=\pi^{\frac
{n+1}{2}}(8(R+1))^{n+1}C^{\bigtriangleup_{\varphi}(B)}[(2e)^{\frac
{n+1}{2}}]$.  Set $D_0=\pi^{n+1}(8(R+1))^{n+1}(C+1)6^{n+1}$.  Clearly $D_0
> D$.  $B \subset N_{m_0}$ is a unital inclusion so by Lemma 3.8
$\bigtriangleup_{\varphi}(B) \leq \bigtriangleup_{\varphi}(N_{m_0}) \leq
\alpha$.  Hence, $n+1-\bigtriangleup_{\varphi}(B) \geq n+1-\alpha$.  
Since $0<\epsilon<1$ \[ \log (\epsilon^{n+1-\bigtriangleup_{\varphi}(B)})+
\log(4^{n+1}D) \leq \log(\epsilon^{n+1-\alpha}) + \log(4^{n+1}D_0). \]
Thus, \begin{eqnarray}\frac {\chi(a_1+\epsilon s_1,\ldots,a_n + \epsilon
s_n, I+ \epsilon s_{n+1}: s_1,\ldots,s_{n+1})}{| \log \epsilon|} & \leq &
-(n+1) + \alpha + \frac {\log(4^{n+1}D_0)}{| \log \epsilon|}.  \nonumber
\end{eqnarray} $D_0$ is independent of $\epsilon$ so by Lemma 3.9
\begin{eqnarray}\delta_0(a_1,\ldots,a_n)=\delta_0(a_1,\ldots,a_n,I) & \leq
& (n+1) + \limsup_{\epsilon \rightarrow 0} \left (-(n+1) + \alpha + \frac
{\log(4^{n+1}D_0)}{| \log \epsilon|} \right) \nonumber \\ & = & \alpha
\nonumber \\ & = & 1 - \sum_{i=1}^{s}
\frac{\alpha_{i}^{2}}{k_{i}^{2}}.\nonumber\end{eqnarray} \end{proof}

\section{Weak Hyperfinite Monotonicity}

	Throughout this section assume $b_1,\ldots,b_p$ are self-adjoint
elements in $M$ and the strongly closed algebra $B$ generated by the $b_j$
is hyperfinite.  We will prove that if $\{b_1,\ldots,b_p \}$ lies in the
$*$-algebra generated by $\{a_1,\ldots,a_n\},$ then
\[\delta_0(b_1,\ldots,b_p) \leq \delta_0(a_1,\ldots,a_n).\]

\noindent This "weak hyperfinite monotonic" inequality has significant
implications in finding sharp lower bounds for $\delta_0(a_1,\ldots,a_n)$
when $M$ is diffuse.

	The argument is simple, despite the notation which shrouds it.   
Because $B$ is hyperfinite matricial microstates of 
$\{b_1,\ldots,b_p \}$ are all approximately unitarily 
equivalent; the proof is nothing more than a trivial generalization of 
Lemma 3.1. It follows that 
$\delta_0(b_1,\ldots,b_p)$ 
reflects the metric entropy of the unitary orbit of a single microstate 
for $\{b_1,\ldots,b_p\}$(provided the microstate 
approximates well enough).  Since the $b_j$ are polynomials of 
the $a_i$ (and thus images of the $a_i$ under Lipschitz maps), the metric 
entropy data carries over to the microstates of 
$\{a_1,\ldots,a_n\}$ and yields lower bounds for the metric entropy of the 
unitary orbit of a microstate for $\{a_1,\ldots,a_n\}$.  Stuffing this 
lower bound information into the modified free entropy dimension machine we 
arrive at the above inequality.    

	In addition to maintaining the conventions set forth in Section 2 
we adopt the following notation in this section:
\begin{itemize}
\item For $r >0$ $(M_{k}^{sa}(\mathbb C))_r$ denotes the operator norm 
ball of $M_{k}^{sa}(\mathbb C)$ of radius $r$ centered at the origin.  
For any $d \in 
\mathbb N$ $((M_{k}^{sa}(\mathbb C))_r)^d$ is the Cartesian product of $d$ 
copies of $(M_{k}^{sa}(\mathbb C))_r$. 
\item For \hspace{.07in}$d\in \mathbb N, K \subset (M_{k}^{sa}(\mathbb C))^d$, and $u 
\in U_k$ define $uKu^*$ to be the set 
\[\{(ux_1u^*,\ldots,ux_du^*):(x_1,\ldots,x_d) \in K\}.\]
\item For $d\in \mathbb N$ and $(x_1,\ldots,x_d)\in 
(M_{k}^{sa}(\mathbb C))^d$ define $U(x_1,\ldots,x_d) = 
\{(ux_1u^*,\ldots,ux_nu^*): u\in U_k\}$. 
\item For $\epsilon > 0, d\in \mathbb N, $ and $S \subset 
(M_{k}^{sa}(\mathbb C))_d$ write $P_{\epsilon}(S)$ for the maximum 
number of points in an $\epsilon$-separated subset of $S$ and 
$\mathcal N_{\epsilon}(S)$ for the $\epsilon$-neighborhood of $S$, both 
taken with respect to the metric $\rho((x_1,\ldots,x_d),(y_1,\ldots,y_d)) = \max 
\{|x_i-y_i|_2 : 1 \leq i \leq d\}.$
\end{itemize}

	First we show that matricial microstates of 
$\{b_1,\ldots,b_p\}$ are approximately unitarily 
equivalent.  We do this with the following two lemmas, the first of which 
makes no use of hyperfiniteness.  

\begin{lemma}If $z_1,\ldots,z_p \in B$ are self-adjoint, 
$\|z_j\| \leq r$ for $1 \leq j\leq p$, and $L, \gamma_0 > 0$, then there 
exist polynomials 
$f_1,\ldots,f_p$ in p noncommuting variables such that for $1 \leq j \leq 
p$:
\begin{itemize}
\item $|f_j(b_1,\ldots,b_p) - z_j|_2 < 2\gamma_0$.
\item $\|f_j(b_1,\ldots,b_p)\| \leq r+1$.
\item For any $k \in \mathbb N$ and $(x_1,\ldots,x_p) \in 
((M_{k}^{sa}(\mathbb C))_{L})^p$ 
$f_j(x_1,...,x_p) \in (M_{k}^{sa}(\mathbb C))_{r+1}.$\end{itemize}
\end{lemma}

\begin{proof}By Kaplansky's Density Theorem there exist polynomials
$g_1,\ldots,g_p$ in $p$ noncommuting variables such that for $1 \leq j
\leq p$: \begin{itemize} \item $|g_j(b_1,\ldots,b_p) - z_j |_2 <
\gamma_0$. \item $\|g_j(b_1,\ldots,b_p)\| \leq r$. \item
$g_j(y_1,\ldots,y_p)$ is self-adjoint for any self-adjoint
operators $y_1,\ldots, y_p.$ \end{itemize}

\noindent There exists an $L_1>L+r$ such that for any $1 \leq j \leq p$ 
and $k \in \mathbb N$ if $(x_1,\ldots,x_p) \in ((M_{k}^{sa}(\mathbb 
C))_L)^p,$
then $\|g_j(x_1,\ldots,x_p)\| \leq L_1$.  Define $f:[-L_1,L_1] \rightarrow
\mathbb R$ by \[f(t) = \left\{ \begin{array}{ll} t & \mbox{if $|t| \leq
r$} \\ r & \mbox{if $r<t \leq L_1$} \\ -r & \mbox{if $-L_1 \leq t < -r$}
\end{array} \right. \] 

\noindent For any $1\leq j \leq p \hspace{.1in}
f(g_j(b_1,\ldots,b_p))=g_j(b_1,\ldots,b_p)$
and: \begin{itemize} \item $|f(g_j(b_1,\ldots,b_p))
- z_j|_2 < \gamma_0$. \item
$\|f(g_j(b_1,\ldots,b_p))\| \leq r$. \item For any
$k \in \mathbb N$ if $(x_1,\ldots,x_p) \in ((M_{k}^{sa}(\mathbb
C))_{L})^p,$ then $f(g_j(x_1,\ldots,x_p))\in (M_{k}^{sa}(\mathbb C))_r$.
\end{itemize} Approximating f uniformly on $[-L_1, L_1]$ by a polynomial 
$h$ (to within sufficiently small $\epsilon > 0$) and setting $f_j=h \circ
g_j$ yields the desired result. \end{proof}

\begin{lemma}If $\epsilon >0$ and $r \geq \max\{\|b_j\|\}_{1\leq j\leq
p}$, then there exist $m \in \mathbb N$ and $\gamma >0$ such that for each
$k \in \mathbb N$ and $(x_1,\ldots,x_p), (y_1,\ldots,y_p) \in
\Gamma_r(b_1,\ldots,b_p;m,k,\gamma)$ there exists a $u\in U_k$ satisfying
\[|ux_ju^*-y_j|_2 < \epsilon\] for $1\leq j \leq p$. \end{lemma}

\begin{proof}By Kaplansky's Density Theorem and the hyperfiniteness of $B$
there exist self-adjoint elements $z_1,\ldots,z_p \in B$ which generate a
finite dimensional algebra and such that $\|z_j\| \leq r$ and $|z_j -
b_j|_2 < \epsilon$ for all $1 \leq j \leq p.$ By the remark following
Corollary 3.4 there exist $m_1 \in \mathbb N$ and $\gamma_1>0$ such that
if $r_0 = \max \{\|z_j\|\}_{1\leq j \leq p}$, $k \in \mathbb N$, and
$(x_1,\ldots,x_p),$ $(y_1,\ldots,y_p)$ $\in \Gamma_{r_0
+1}(z_1,\ldots,z_p;m_1,k,\gamma_1)$, then there exists a $u \in U_k$
satisfying $|ux_ju^*-y_j|_2 < \epsilon$ for $1 \leq j \leq p$.

	By Lemma 4.1 for $\frac {\epsilon}{2} > \gamma_0 > 0$ there exist
polynomial $f_1,\ldots,f_p$ in $p$ noncommuting variables such that for $1
\leq j \leq p$:  
\begin{itemize} \item $|f_j(b_1,\ldots,b_p) - z_j|_2 <
2\gamma_0 < \epsilon$. 
\item $\|f_j(b_1,\ldots,b_p)\| \leq r_0+1$. 
\item For any $k \in \mathbb N$ if $(x_1,\ldots,x_p) \in 
((M_{k}^{sa}(\mathbb
C))_r)^p,$ then $f_j(x_1,\ldots,x_p) \in (M_{k}^{sa}(\mathbb C))_{r_0+1}.$
\end{itemize} By making $\gamma_0$ sufficiently small it follows that for
any $1 \leq j \leq m_1$ and $1 \leq i_1,\ldots,i_j \leq p$ 

\[|f_{i_1}(b_1,\ldots,b_p) \cdots
f_{i_j}(b_1,\ldots,b_p) - z_{i_1} \cdots z_{i_j} |_2
< \gamma_1.\] 

\noindent Hence by choosing $m \in \mathbb N$ large enough
and $\gamma >0$ small enough if $k \in \mathbb N$ and
$(x_1,\ldots,x_p) \in 
\Gamma_r(b_1,\ldots,b_p;m,k,\gamma)$,
then $|f_j(x_1,\ldots,x_p)-x_j|_2 \leq
|f_j(b_1,\ldots,b_p) - b_j|_2 + \epsilon$
for $1 \leq j \leq p$ and
\[(f_1(x_1,\ldots,x_p),\ldots,f_p(x_1,\ldots,x_p)) \in
\Gamma_{r_0+1}(z_1,\ldots,z_p;m_1,k,\gamma_1).\]

	Finally suppose $k \in \mathbb N$ and $(x_1,\ldots,x_p), 
(y_1,\ldots,y_p) \in 
\Gamma_r(b_1,\ldots,b_p;m,k,\gamma).$  
For any $1 \leq j \leq p$

\[|f_j(x_1,\ldots,x_p) - x_j|_2 \hspace{.07in} 
\leq \hspace{.07in} 
|f_j(b_1,\ldots,b_p) 
- b_j|_2 + \epsilon \hspace{.07in} \leq \hspace{.07in} |z_j - 
b_j|_2 + 2\epsilon \hspace{.07in}
< \hspace{.07in}3\epsilon.\]
Similarly for $1 \leq j \leq p \hspace{.1in} |f_j(y_1,\ldots,y_p) - y_j|_2 
< 3\epsilon$.  By the preceding two paragraph there exists a $u \in U_k$ 
such that for $1 \leq j \leq p \hspace{.1in} |uf_j(x_1,\ldots,x_p)u^* - 
f_j(y_1,\ldots,y_p)|_2 < 
\epsilon$.  So for $1 \leq j \leq p \hspace{.1in} |ux_ju^* - y_j|_2 < 
7\epsilon$. \end{proof}

	In the next lemma suppose $ r > \max \{ \|b_j\| \}_{1 
\leq j \leq p}$.

\begin{lemma}For each $0< \epsilon <1$ there exist corresponding 
$m_{\epsilon} \in  
\mathbb N$ and $\gamma_{\epsilon} > 0$ such that if $L>0$ and 
$\langle (x_1^{(k)},\ldots,x_p^{(k)}) \rangle_{k=1}^{\infty}$ is a 
sequence 
satisfying 
$(x_1^{(k)},\ldots,x_p^{(k)}) \in (M_{k}^{sa}(\mathbb C))^p$ for all k 
and $(x_1^{(k)},\ldots,x_p^{(k)}) \in 
\Gamma_r(b_1,\ldots,b_p;m_{\epsilon},k,\gamma_{\epsilon})$ 
for sufficiently large k, then 
\[\limsup_{k \rightarrow \infty}[k^{-2} \cdot \log(P_{\epsilon 
L}(U(x_{1}^{(k)},\ldots,x_p^{(k)})))] \geq 
\chi_r(b_1+\epsilon 
s_1,\ldots,b_p+\epsilon s_p: s_1,\ldots,s_p) + p | \log 
\epsilon| - K_0\]where $K_0 = p \cdot \log((2+L)\sqrt{2\pi e})$.
\end{lemma} 

\begin{proof}Suppose $0< \epsilon <1$.  By Lemma 4.2 there exist 
$m_{\epsilon} \in  \mathbb N$ and $\gamma_{\epsilon} > 0$ such that if $k 
\in \mathbb N,$ 

\[(z_1,\ldots,z_p) \in \Gamma_r(b_1 +\epsilon s_1,\ldots,b_p +\epsilon 
s_p:s_1,\ldots,s_p;m_{\epsilon},k,\gamma_{\epsilon}),\] 

\noindent and  $(x_1,\ldots,x_p) \in
\Gamma_r(b_1,\ldots,b_p;m_{\epsilon},k,\gamma_{\epsilon})$, then there
exists a $u\in U_k$ satisfying $|ux_ju^*-z_j|_2 < 2\epsilon$ for $1 \leq j
\leq p$.

	Assume $\langle (x_1^{(k)},\ldots,x_p^{(k)})  
\rangle_{k=1}^{\infty}$ satisfies the hypothesis of the lemma with
$m_{\epsilon}$ and $\gamma_{\epsilon}$ chosen according to the preceding
paragraph.  For sufficiently large $k$

\[\Gamma_r(b_1 +\epsilon s_1,\ldots,b_p+\epsilon
s_p:s_1,\ldots,s_p;m_{\epsilon},k,\gamma_{\epsilon}) \subset 
\mathcal N_{2\epsilon}(U(x_1^{(k)},\ldots,x_p^{(k)})).\]
Find an $\epsilon L$-separated set $W_k$ of 
$U(x_1^{(k)},\ldots,x_p^{(k)})$ 
(with respect to the $\rho$ metric) 
of maximum cardinality.  
$\mathcal N_{2\epsilon}(U(x_{1}^{(k)},\ldots,x_p^{(k)}))\subset 
\mathcal N_{(2+L)\epsilon}(W_k)$.  For large enough $k$

\begin{eqnarray} \text{vol}(\Gamma_r(b_1 +\epsilon 
s_1,\ldots,b_p+\epsilon  
s_p:s_1,\ldots,s_p;m_{\epsilon},k,\gamma_{\epsilon})) & \leq & 
\text{vol}(\mathcal N_{(2+L)\epsilon}(W_k)) \nonumber \\ & \leq & |W_k| 
\cdot 
\frac {\pi^{\frac {pk^2}{2}}((2+L)\epsilon\sqrt{k})^{pk^2}}{\Gamma 
\left (\frac {k^2}{2} + 1 \right)^p}. \nonumber
\end{eqnarray}
	
	By the preceding inequality $\chi_r(b_1 +\epsilon
s_1,\ldots,b_p+\epsilon
s_p:s_1,\ldots,s_p;m_{\epsilon},\gamma_{\epsilon})$ is dominated by

\begin{eqnarray*} & & \limsup_{k \rightarrow \infty} \left[k^{-2} \cdot
\log|W_k| + p \log((2+L)\epsilon\sqrt{\pi k}) - p k^{-2} \cdot \log
\left (\Gamma \left (\frac {k^2}{2} +1\right) \right)  + \frac {p}{2}
\log k \right] \\ & \leq & \limsup_{k \rightarrow \infty}
\left [ k^{-2} \cdot \log|W_k| +p \log((2+L)\epsilon\sqrt{\pi }) -
p k^{-2} \log \left (\left ( \frac {k^2}{2e} \right)^{\frac {k^2}{2}}
\right) + p \log k \right] \\ & = & \limsup_{k \rightarrow
\infty} \left [k^{-2} \cdot \log|W_k| +p \log((2+L)\epsilon\sqrt{\pi
}) - \frac {p}{2} \cdot \log \left (\frac {k^2}{2e} \right) + p 
\log k \right] \\ & = & \limsup_{k\rightarrow \infty}[k^{-2}
\cdot \log|W_k| +p \log((2+L)\epsilon\sqrt{\pi })+ \frac {p}{2}
\log (2e)] \\ & = & p \log((2+L)\sqrt{2\pi e}) + 
\log \epsilon + \limsup_{k \rightarrow \infty} (k^{-2} \cdot
\log|W_k|) \\ & = & K_0 + p \log \epsilon + \limsup_{k
\rightarrow \infty} [k^{-2} \cdot \log(P_{\epsilon
L}(U(x_{1}^{(k)},\ldots,x_p^{(k)})))]. \\ \end{eqnarray*}

By the above calculation $\chi_r(b_1 +\epsilon 
s_1,\ldots,b_p+ 
\epsilon s_p:s_1,\ldots,s_p) + p |\log \epsilon| - K_0$ is dominated 
by \[\limsup_{k \rightarrow \infty}[k^{-2} \cdot \log(P_{\epsilon
L}(U(x_1^{(k)},\ldots,x_p^{(k)})))].
\]
\end{proof}

	Write $A$ for the $*$-algebra generated by $\{a_1,\ldots,a_n\}$.

\begin{lemma}If  $\{b_1,\ldots,b_p\} \subset A$, then 
there exists an $L>0$ such that for any $0 < \epsilon < 1, m \in \mathbb N$, and 
$\gamma > 0$ there is a sequence 
$\langle (x_{1}^{(k)},\ldots,x_{n}^{(k)}) \rangle_{k=1}^{\infty}$ 
satisfying $(x_{1}^{(k)},\ldots,x_{n}^{(k)}) \in (M_{k}^{sa}(\mathbb C))^n$ for all 
$k$, $(x_{1}^{(k)},\ldots,x_{n}^{(k)}) \in 
\Gamma_{R+1}(a_1,\ldots,a_n;m, 
k, \gamma)$ for sufficiently large $k$, and
\[ \limsup_{k \rightarrow \infty} [k^{-2} \cdot 
\log(P_{4\epsilon\sqrt{n}}(U(x_{1}^{(k)},\ldots,x_{n}^{(k)})))] \geq   
\chi_{\lambda}(b_1 + \epsilon s_1,\ldots,b_p + 
\epsilon s_p: s_1,\ldots,s_p) + p |\log \epsilon| - K_1 \]
where $K_1=p \cdot \log((2+4\sqrt{n}L)\sqrt{2\pi e})$ $\hspace{.01in}$ and 
$\hspace{0.01in}$  $\lambda= L(R+1)+ \max\{\|b_j \|\}_{1 \leq j 
\leq p}$.
\end{lemma} 

\begin{proof}There exist polynomials $f_1,\ldots,f_p$ in $n$ noncommuting
variables and with no constant terms such that for $1 \leq j \leq p
\hspace{.06in} f_j(a_1,\ldots,a_n) = b_{j}$ and such that $f_j$ of an
$n$-tuple of self-adjoint operators is a self-adjoint element.  There
exists a constant $L>0$ such that if $k\in \mathbb N$ and
$\xi_1,\ldots,\xi_n,\eta_1,\ldots,\eta_n \in (M_{k}^{sa}(\mathbb 
C))_{R+1}$, then
for all $1 \leq j \leq p$ \[|f_j(\xi_1,\ldots,\xi_n) - 
f_j(\eta_1,\ldots,\eta_n)|_2
\leq L \cdot \max \{|\xi_i - \eta_i|_2: 1 \leq i \leq n \} \]

\noindent and $\|f_j(\xi_1,\ldots,\xi_n)\| \leq L(R+1)$.

	Suppose $0<\epsilon<1, m \in \mathbb N$, and $\gamma > 0.$ By
Lemma 4.3 there exist an $m_{\epsilon} \in \mathbb N$ and
$\gamma_{\epsilon} > 0$ such that if $\langle
(x_1^{(k)},\ldots,x_p^{(k)}) \rangle_{k=1}^{\infty}$ is a sequence
satisfying $(x_1^{(k)},\ldots,x_p^{(k)}) \in (M_{k}^{sa}(\mathbb C))^p$
for all $k$ and $(x_1^{(k)},\ldots,x_p^{(k)})\in
\Gamma_{\lambda}(b_1,\ldots,b_p;m_{\epsilon},k,\gamma_{\epsilon})$ for
sufficiently large $k,$ then

\[ \limsup_{k \rightarrow \infty}[k^{-2} \cdot \log(P_{4 
\epsilon \sqrt{n} L}(U(x_{1}^{(k)},\ldots,x_p^{(k)})))] \geq  
\chi_{\lambda}(b_1+\epsilon   
s_1,\ldots,b_p+\epsilon s_p: s_1,\ldots,s_p) + p|\log 
\epsilon| -K_1.\]

\noindent where $K_1=p \cdot \log((2+4\sqrt{n} L)\sqrt{2\pi e})$.  By the 
assumed existence of finite dimensional approximants for 
$\{a_1,\ldots,a_n\}$ there exists a 
$k_0 \in \mathbb N$ such that for each $k \geq k_0$ there is an 
$(x_{1}^{(k)},\ldots,x_{n}^{(k)}) \in \Gamma_{R+1}(a_1,\ldots,a_n; 
m,k,\gamma)$ satisfying 

\[(f_1(x_1^{(k)},\ldots,x_n^{(k)}),\ldots,f_p(x_1^{(k)},\ldots,x_n^{(k)})) 
\in 
\Gamma_{\lambda}(b_1,\ldots,b_p;m_{\epsilon},k,\gamma_{\epsilon}).\]
For each $k \geq k_0$ and $1 \leq j \leq p$ set 
$y_{j}^{(k)}=f_j(x_1^{(k)},\ldots,x_n^{(k)})$.  It follows that 

\[ \limsup_{k \rightarrow \infty}[k^{-2}\cdot
\log(P_{4\epsilon\sqrt{n} L}(U(y_{1}^{(k)},\ldots,y_p^{(k)})))] \geq 
\chi_{\lambda}(b_1+\epsilon s_1,\ldots,b_p +\epsilon
s_p: s_1,\ldots,s_p) + p| \log \epsilon| -K_1. \]

	For $k \geq k_0$ and any $u, v \in U_k$, $1 \leq j \leq p,$ 
observe that $L \cdot \max \{|ux_i^{(k)}u^* - vx_i^{(k)}v^*|_2 : 1 \leq
i \leq n \}$ dominates

\[ |q_j(ux_1^{(k)}u^*, \ldots, ux_n^{(k)}u^*) 
-q_j(vx_1^{(k)}v^*,\ldots,vx_n^{(k)}v^*)|_2  
|uy_j^{(k)}u^* -vy_j^{(k)}v^*|_2. 
\]
It follows that for any $k \geq k_0$ 
\[k^{-2} \cdot \log (P_{4 \epsilon \sqrt{n}} (U(x_1^{(k)}, \ldots, 
x_n^{(k)}))) \geq k^{-2} \cdot \log (P_{4 \epsilon \sqrt{n} 
L}(U(y_1^{(k)}, \ldots, y_p^{(k)}))). \]

\noindent By the last sentence of the preceding paragraph we're done.
\end{proof} 

\begin{theorem}(Weak Hyperfinite Monotonicity) If 
$\{b_1,\ldots,b_p\} \subset A$, 
then

\[\delta_0(b_1,\ldots,b_p) \leq 
\delta_0(a_1,\ldots,a_n).\]
\end{theorem}

\begin{proof}Consider the constants $L$ and $\lambda$ corresponding to 
$b_1,\ldots,b_p$ in Lemma 4.4.  Suppose $0< \gamma, 
\epsilon < \frac {1}{2}$ and $m \in \mathbb N.$  By Lemma 4.4 
there exists a sequence 
$<(x_{1}^{(k)},\ldots,x_{n}^{(k)})>_{k=1}^{\infty}$ 
satisfying $(x_{1}^{(k)},\ldots,x_{n}^{(k)}) \in (M_{k}^{sa}(\mathbb 
C))^n$ for all $k\in \mathbb N$, $(x_{1}^{(k)},\ldots,x_{n}^{(k)}) \in 
\Gamma_{R+1}\left(a_1,\ldots,a_n;m,k,\frac 
{\gamma}{(8(R+2))^{m}} \right)$ for 
sufficiently large $k$, and 

\[\limsup_{k \rightarrow \infty} [k^{-2} \cdot
\log(P_{4\epsilon\sqrt{n}}(U(x_{1}^{(k)},\ldots,x_{n}^{(k)})))] \geq
\chi_{\lambda}(b_1 + \epsilon s_1,\ldots,b_p + 
\epsilon s_p: s_1,\ldots,s_p) + p | \log \epsilon| - 
K_1 \]where $K_1 = p \cdot \log((2+4\sqrt{n}L)\sqrt{2\pi e})$.

	By Corollary 2.14 of [11] there is an $N \in \mathbb N$ such that
if $k \geq N$ and $\sigma$ is a Radon probability measure on
$((M_{k}^{sa}(\mathbb C))_{R+1})^{2n}$ invariant under the action

\[(\xi_1,\ldots,\xi_n,\eta_1,\ldots,\eta_n) \mapsto
(\xi_1,\ldots,\xi_n,u\eta_1 u^{*},\ldots,u\eta_n u^{*})\] 

\noindent for $u\in U_k$, then
$\sigma(\omega_k) > \frac {1}{2}$ where 

\begin{eqnarray}
\omega_k=\{(\xi_1,\ldots,\xi_n,\eta_1,\ldots,\eta_n)\in
((M_{k}^{sa}(\mathbb C))_{R+1})^{2n} & : & \{\xi_1,\ldots,\xi_n\} \text {
and } \{\eta_1,\ldots,\eta_n\} \nonumber \\ & & \text { are } \left( m,
\frac {\gamma}{4^m}\right) -\text {free}\}. \nonumber \end{eqnarray} For
$k\in \mathbb N$ write $\nu_k$ for the atomic probability measure
concentrated at $(x_1^{(k)},\ldots,x_n^{(k)})$ and $m_k$ for the
probability measure obtained by restricting vol to
$\Gamma_{2\epsilon}\left (\epsilon s_1,\ldots,\epsilon s_n;  m,k,\frac
{\gamma}{8^m}\right)$ and normalizing appropriately.  $\nu_k \times m_k$
is a Radon probability measure on $((M_{k}^{sa}(\mathbb C))_{R+1})^{2n}$
invariant under the $U_k$-action described above.  Write $F_k$ for the set
of all $(z_1,\ldots,z_n) \in \Gamma_{2\epsilon} \left (\epsilon
s_1,\ldots,\epsilon
s_n;m,k,\frac{\gamma}{8^m}\right)$ such that
$\{z_1,\ldots,z_n\}$ and $\{x_{1}^{(k)},\ldots,x_{n}^{(k)}\}$ are $\left(
m,\frac {\gamma}{4^m} \right )$-free.

	For $k \geq N$ $\frac {1}{2} < (\nu_k \times m_k)(\omega_k) =
m_k(F_k)$.  Set $E_k = (x_{1}^{(k)},\ldots,x_{n}^{(k)})  + F_k$.  For $u
\in U_k,$ $\text{vol}(uE_ku^*) = \text{vol}(E_k)$ and $uE_ku^*$ is
contained in \[ \Gamma_{R+1+2\epsilon}(a_1 + \epsilon s_1,...,a_n +
\epsilon s_n: \epsilon s_1,\ldots,\epsilon s_n; m, k, \gamma).\] For
each $k\in N$ there exists a subset $\langle u_{k,s}\rangle_{s \in S_k}$
of $U_k$ such that $|S_k| =
P_{4\epsilon\sqrt{n}}(U(x_{1}^{(k)},\ldots,x_{n}^{(k)}))$ and for $s,
s^{\prime} \in S_k$ with $s \neq s^{\prime}$

\[ \max\{|u_{k,s}x_{i}^{(k)}u_{k,s}^{*} - 
u_{k,s^{\prime}}x_{i}^{(k)}u_{k,s^{\prime}}^{*}|_2 : 1 \leq i \leq n \} > 
4\epsilon\sqrt{n}. \]

\noindent Since $F_k$ is $\| \cdot \|_2$-bounded by $2\epsilon\sqrt{nk}
\hspace{.1in} (u_{k,s}E_ku_{k,s}^*)  \bigcap
(u_{k,s^{\prime}}E_ku_{k,s^{\prime}}^{*}) = \emptyset$ for $s, s^{\prime}
\in S_k, s \neq s^{\prime}$.  Hence for $k \geq \ N \hspace{.05in}
\text{vol}(\Gamma_{R+1+2\epsilon}(a_1 + \epsilon s_1,...,a_n + \epsilon 
s_n:\epsilon s_1,\ldots, \epsilon s_n; m, k, \gamma))$
dominates

\begin{eqnarray} \text{vol}(\bigsqcup_{s\in S_k} u_{k,s}E_ku_{k,s}^*) = 
|S_k| \cdot \text{vol}(F_k) & = & |S_k| \cdot m_k(F_k) \cdot 
\text{vol} (\Gamma_{2\epsilon} (\epsilon s_1,\ldots,\epsilon s_n; 
m, k, \gamma /8^m)) 
\nonumber \\ & > & 
1/2 \cdot |S_k| \cdot 
vol(\Gamma_{2\epsilon}(\epsilon s_1,\ldots,\epsilon s_n; 
m, k, \gamma / 8^m)).\nonumber\end{eqnarray}  

	By the last sentence of the preceding paragraph \begin{eqnarray*}
& & \chi_{R+1+2\epsilon}(a_1 + \epsilon s_1,...,a_n + \epsilon s_n :
\epsilon s_1,\ldots, \epsilon s_n; m, \gamma) \\ & \geq & \limsup_{k
\rightarrow \infty} (k^{-2} \cdot [\log(1/2 \cdot |S_k| \cdot \text{vol}
(\Gamma_{2\epsilon}(\epsilon s_1,\ldots,\epsilon s_n;  m, k, \gamma/
8^m)))] + n/2 \cdot \log k)  \\ & \geq & \limsup_{k \rightarrow \infty}
(k^{-2} \cdot \log(|S_k|)) + \liminf_{k \rightarrow \infty} (k^{-2}
\cdot [\log(\text{vol} (\Gamma_{2\epsilon}(\epsilon s_1,...,\epsilon s_n; 
m, k,\gamma/ 8^m)))] + n/2 \cdot \log k)  \\ & \geq & \limsup_{k \rightarrow
\infty} (k^{-2} \cdot \log(|S_k|)) + \chi_{2\epsilon}(\epsilon
s_1,\ldots,\epsilon s_n) \\ & = & \limsup_{k \rightarrow \infty} [k^{-2}
\cdot \log(P_{4\epsilon \sqrt {n}}(U(x_{1}^{(k)},\ldots,x_{n}^{(k)})))] +
n \log(\epsilon \sqrt {2\pi e} ) \\ & \geq & \chi(b_1 + \epsilon
s_1,\ldots,b_p + \epsilon s_p: s_1,\ldots,s_p) + p | \log \epsilon |
+ n \log(\epsilon \sqrt{2\pi e}) - K_1.  \end{eqnarray*} where we
used regularity of $\{\epsilon s_1,\ldots,\epsilon s_n\}$ going from the
third to the fourth lines above.  $m$ and $\gamma$ being arbitrary it
follows that 

\begin{eqnarray*} \chi(a_1 + \epsilon s_1,\ldots,a_n +
\epsilon s_n: s_1,\ldots,s_n)  & \geq & \chi(b_1+ \epsilon
s_1,...,b_p+ \epsilon s_p:  s_1,\ldots,s_p) \\ & & +
\hspace{.05in}(p - n) \cdot | \log \epsilon | + \hspace {.05in} n \cdot
\log (\sqrt {2\pi e}) - K_1. \end{eqnarray*} Dividing by $|\log
\epsilon|,$ taking $\limsup$'s as $\epsilon \rightarrow 0$, and adding $n$
to both sides yields \[ \delta_0(a_1,\ldots,a_n) \geq
\delta_0(b_1,\ldots,b_p).\] \end{proof}

\begin{corollary}If $a\in M$, then $\delta_0(a_1,\ldots,a_n) \geq
\delta_0(a).$ \end{corollary} \begin{proof}Find a sequence $\langle z_k
\rangle_{k=1}^{\infty}$ in $A$ such that $z_k \rightarrow a$ strongly.  
By Proposition 6.14 of [9] and Corollary 6.7 of [10] $\liminf_{k
\rightarrow \infty} \delta_0(z_k) = \liminf_{k \rightarrow \infty}
\delta(z_k) \geq \delta (a) = \delta_0(a)$.  For each $k \hspace{.1in}
z_k$ generates a hyperfinite von Neumann algebra; by Lemma 4.2 for each
$k \hspace{.07in}\delta_0(a_1,\ldots,a_n) \geq \delta_0(z_k)$ so the
preceding sentence yields the desired result.  \end{proof}
\begin{corollary}If $M$ has a diffuse von Neumann subalgebra, then
$\delta_0(a_1,\ldots,a_n) \geq 1$. \end{corollary} \begin{proof}Find a
maximal abelian subalgebra $N$ of the diffuse von Neumann subalgebra.  
$N$ has a self-adjoint generator $a$.  $N$ must be diffuse since it is a
maximal abelian subalgebra of a diffuse von Neumann algebra.  
Consequently $a$ has no eigenvalues.  Apply Corollary 4.6.  \end{proof}

\begin{remark} By [10] if $M$ has a regular diffuse von Neumann
subalgebra, then $\delta_0(a_1,\ldots,a_n) \leq 1.$ By [4] if there
exists a sequence of Haar unitaries $\langle u_j \rangle_{j=1}^s$ such
that the sequence generates $M$ as a von Neumann algebra and for each $j
\in \mathbb N$ $u_{j+1}u_j u_{j+1}^* \in \{u_1,\ldots,u_j\}^{\prime
\prime},$ then $\delta_0(a_1,\ldots,a_n) \leq 1.$ Combining these
results with Corollary 4.7, it follows that for any self-adjoint
generators $a_1,\ldots,a_n$ for an $M$ which satisfies either of the two
conditions and which is also embeddable into the ultraproduct of the
hyperfinite $II_1$-factor, $\delta_0(a_1,\ldots,a_n)=1.$ In other words,
$\delta_0$ is a von Neumann algebra invariant for such algebras.  In
particular, $\delta_0(M) =1$ when $M$ can be embedded into the
ultraproduct of the hyperfinite $II_1$-factor and $M$ has a Cartan
subalgebra, $M = N_1 \otimes N_2$ for $II_1$-factors $N_1$ and $N_2,$
or when $M$ is a group von Neumann algebra associated to the groups
$SL_n(\mathbb Z), n \geq 3.$ \end{remark}

\section{Lower Bound for Finite Dimensional Algebras}

In this section we calculate the lower bound for
$\delta_0(a_1,\ldots,a_n)$ when $M$ is finite dimensional.  Without loss
of generality assume throughout this section that $M = \oplus_{i=1}^{p}
M_{k_i}(\mathbb C)$ and $\varphi = \oplus_{i=1}^{p} \alpha_i tr_{k_i}$
where $p \in \mathbb N$ and $\alpha_i >0$ for each $i.$ The first lemma we
present is not necessary but it's convenient.

\begin{lemma}There exists an $x \in M$ such that the $*$-algebra 
generated by $x$ is $M$.\end{lemma}  

	The proof is not hard and we omit it.

	As in the preceding section the calculation of the lower bound 
amounts to looking at the packing number of unitary orbits of microstates.  
We use two ingredients.

	For a representation $\pi: M \rightarrow M_k(\mathbb C)$ define 
$H_{\pi}$ to be the unitary group of $(\pi(M))^{\prime}$ and 
$X_{\pi}=U_k/H_{\pi}$.  Endow $X_{\pi}$ with the quotient metric from 
the $|\cdot|_2$-metric on $U_k$.  Call this metric on $X_{\pi} 
\hspace {.1in} d_{\pi}$.  The first ingredient is a packing number 
estimate for certain homogeneous spaces $X_{\pi}$.

\begin{lemma}There exists a $\kappa > 0$ with the property that for every
$\varepsilon > 0$ there is a corresponding sequence
$\langle \sigma_k \rangle_{k=1}^{\infty}$ such that for each $k$ 
$\sigma_k : M
\rightarrow M_k(\mathbb C)$ is a $*$-homomorphism and for $k$ sufficiently
large: \begin{itemize} \item $\|tr_k \circ \sigma_k - \varphi \| <
\varepsilon.$ \item For each $k$ setting $H_k = H_{\sigma_k}$ and $X_k =
X_{\sigma_k}$ we have that $H_k$ is a tractable Lie subgroup of $U_k$
satisfying $k^2(\bigtriangleup_{\varphi}(M) - \varepsilon) \leq
\dim(X_k)$.  \item For any $\epsilon > 0$

\[ \left(\frac{\kappa}{\epsilon}\right)^{\dim X_k} \leq 
P(X_k, \epsilon) \]

\noindent where $P(X_k, \epsilon)$ is the maximum number of points in an 
$\epsilon$-separated subset of $X_k.$
\end{itemize}
\end{lemma} 

	We quarantine the proof of Lemma 5.2 to the Addendum, merely
noting for now that the argument will require some technical modifications
to the proofs in [7].

	From now on fix $x$ as in Lemma 5.1.  Given a representation 
$\pi:M \rightarrow M_k(\mathbb C)$ define $U_{\pi}(x) = \{u\pi(x)u^*: u\in U_k\}$ 
and endow $U_{\pi}(x)$ with the inherited $| \cdot |_2$-metric.  For $u 
\in U_k$ denote $\dot{u}$ to be the image of $u$ in $X_{\pi}$ and define 
$f_{\pi}: U_{\pi}(x) \rightarrow X_{\pi}$ by $f_{\pi}(u\pi(x)u^*) = \dot{u}$.  $f_{\pi}$ is 
well-defined for if $u,v\in U_k, \dot{u} = \dot{v} \Longleftrightarrow 
v^*u \in H_{\pi} \Longleftrightarrow u\pi(x)u*=v\pi(x)v*$.

	For the second ingredient recall that in Section 3 covering number 
estimates with respect to the induced operator norm metrics yield the 
desired upper bounds for $\delta_0(a_1,\ldots,a_n)$.  Part of the 
explanation for this lies in the trivial observation that if $u,v\in U_k$ 
and $z\in M_k(\mathbb C)$, then $|uzu^*-vzv^*|_2 \leq 2 \|u-v\| \cdot 
|z|_2$.  The second ingredient more or less says the reverse: there exists 
a constant $L > 0$ such that $d_{\pi}(\dot{u},\dot{v}) \leq L \cdot |u\pi(x)u^* - 
v\pi(x)v^*|_2$.

\begin{lemma}If $z,p\in M_k(\mathbb C)$ with $p$ a projection and $zz^*, 
z^*z \leq \|z^*z\|p$, then there exists a $y \in M_k(\mathbb C)$ 
satisfying $yy^*=y^*y=p$ and
\[ |y-z|_2 \leq |p-z^*z|_2 + |p-e|_2 \leq 2 |p-z^*z|_2 \]
where e is the projection onto the range of $z^*z$.
\end{lemma}
\begin{proof}Denote the polar decomposition of $z$ by $z=u|z|$ and use 
the spectral theorem to write $|z|= \sum_{j=1}^{m} \beta_j e_j$ where the 
$e_j$ are mutually orthogonal rank one projections satisfying $e_1+\cdots+e_m=p$ and $\beta_j \geq 0$.  Now estimate:

\begin{eqnarray}|z-u|_{2}^{2} = |u|z| - up|_{2}^{2} \leq ||z|-p|_{2}^{2} 
= \frac{1}{k} \cdot \sum_{j=1}^{m}(1-\beta_j)^2 & \leq & 
\frac{1}{k} \cdot \sum_{j=1}^{m} (1-\beta_j)^2(1+\beta_j)^2 \nonumber \\ & 
= & |p-z^*z|_{2}^{2}. \nonumber \end{eqnarray}
$|z-u|_2 \leq |p-z^*z|_2$.  Since $uu^*, u^*u \leq p$ there 
exists a partial isometry $v$ such that $vv^*=p-uu^*$ and $v^*v=p-u^*u$.  
So if $y=u+v$, then $yy^*=y^*y=p$.  $u^*u=e$ whence

\begin{eqnarray}|z-y|_2 \leq |z-u|_2 + |v|_2 \leq |p-z^*z|_2 + 
(tr_k(v^*v))^{ \frac{1}{2} } & = & |p-z^*z|_2 + 
(tr_k(p-u^*u))^{ \frac{1}{2} } 
\nonumber \\ & = & |p-z^*z|_2 + |p-e|_2 \nonumber \\ & \leq & 2|p-z^*z|_2. 
\nonumber 
\end{eqnarray}                                                                                                                                                                         
\end{proof}
 	Using Lemma 5.3 we obtain the second ingredient:
\begin{lemma}$\{f_{\pi}: \text {for some} \hspace{.1in} k\in \mathbb N 
\hspace{.1in} \pi:M\rightarrow 
M_k(\mathbb C) \hspace{.1in}  \text{is a representation} \}$ is uniformly 
Lipschitz.
\end{lemma}
\begin{proof}Suppose $\pi:M \rightarrow M_k(\mathbb C)$ is a 
representation.  Because $| \cdot|_2$ is unitarily invariant it suffices 
to show that for any $u \in U_k$
\[\inf_{h\in H_{\pi}} |u-h|_2 = d_{\pi}(f_{\pi}(u\pi(x)u^*), 
f_{\pi}(\pi(x))) \leq L \cdot |u\pi(x)u^* - \pi(x)|_2\]
where $L>0$ is a constant dependent only on $x$.

	If $u \in U_k$ then set $\epsilon = |u\pi(x)u^* - \pi(x)|_2$.  
Denote $<e_{jl}^{(i)}>_{1\leq i\leq p, 1\leq j,l \leq k_i}$ to be the 
canonical s.m.u. for $M$.  There exist polynomials in two noncommuting 
variables $<q_{jl}^{(i)}>_{1\leq i\leq p, 1 \leq 
j,l \leq k_i}$ such that for any $i,j,$ and $l$ $q_{jl}^{(i)}(x,x^*)=e_{jl}^{(i)}$.  
Set $y_{jl}^{(i)}=\pi (e_{jl}^{(i)})$.  There exists a constant $C > 0$ 
dependent only on $x$ such that for any $i,j,$ and $l$

\[|uy_{jl}^{(i)} - y_{jl}^{(i)}u|_2 = |uy_{jl}^{(i)}u^* - y_{jl}^{(i)}|_2 
\leq C\epsilon.\]
Set $K= \sum_{i=1}^{p} k_i$.  By the above inequality $|u\pi(I)u^* - 
\pi(I)|_2 < CK\epsilon$.  Setting $f$ to be the projection 
onto the orthogonal complement of the range of $\pi(I), |ufu^*-f|_2 < 
CK\epsilon$.  Now

\begin{eqnarray*} |u - [(\sum_{1\leq i \leq p, 1\leq j \leq k_i}
y_{jj}^{(i)}uy_{jj}^{(i)}) + fuf]|_2 & \leq & (\sum_{1\leq i
\leq p, 1\leq j \leq k_i} |uy_{jj}^{(i)} - y_{jj}^{(i)}u|_2 \cdot
\|y_{jj}^{(i)}\|) + |uf-fu|_2 \cdot \|f\| \\ & \leq & 2CK\epsilon.
\end{eqnarray*}

\noindent For any $1\leq i \leq p, 1\leq j,l \leq k_i$,
\begin{eqnarray}|y_{jj}^{(i)}uy_{jj}^{(i)}-y_{jl}^{(i)}uy_{lj}^{(i)}|_2 & 
\leq & |y_{jj}^{(i)}uy_{jj}^{(i)}-y_{jj}^{(i)}u|_2 + 
|y_{jl}^{(i)}y_{lj}^{(i)}u - y_{jl}^{(i)}uy_{lj}^{(i)}|_2 \nonumber \\ & 
\leq & |uy_{jj}^{(i)}-y_{jj}^{(i)}u|_2 + |y_{lj}^{(i)}u - uy_{lj}^{(i)}|_2 
\nonumber \\ & \leq & 2C\epsilon. \nonumber
\end{eqnarray}
By Lemma 5.3 there exists for each $1\leq i \leq p$ a $v_i \in M_k(\mathbb 
C)$ such that $v_iv_{i}^{*}=v_{i}^{*}v_{i}=y_{11}^{(i)}$ and $|v_i - 
y_{11}^{(i)}uy_{11}^{(i)}|_2 \leq 
2|y_{11}^{(i)}uy_{11}^{(i)}u^*y_{11}^{(i)} - y_{11}^{(i)}|_2$.  So 

\[ |v_i-y_{11}^{(i)}uy_{11}^{(i)}|_2 \leq 2 
\|y_{11}^{(i)}\|^2 \cdot |uy_{11}^{(i)}u^* - y_{11}^{(i)}|_2 \leq 
2C\epsilon.\] 
Similarly there exists a $v \in M_k(\mathbb C)$ such that $v^*v=vv^*=f$ 
and $|v-fuf|_2 \leq 2|fufu^*f-f|_2 \leq 2 CK
\epsilon$.	

	Consider $z=(\sum_{1 \leq i \leq p, 1 \leq j \leq k_i} 
y_{j1}^{(i)}v_i y_{1j}^{(i)}) + v$.  It's easy to check that $z$ is a 
unitary and because $z$ commutes with all the $y_{jl}^{(i)}, z\in 
H_{\pi}$.  Finally by the last three inequalities of the preceding 
paragraph,

\begin{eqnarray*}|u-z|_2 & \leq & | u - [(\sum_{1\leq i \leq p, 1\leq j 
\leq k_i}y_{jj}^{(i)}uy_{jj}^{(i)}) + fuf]|_2 + \left (\sum_{1\leq i \leq 
p, 1\leq j \leq k_i}|y_{jj}^{(i)}uy_{jj}^{(i)} - 
y_{j1}^{(i)}v_iy_{1j}^{(i)}|_2 \right) \\
& & \\ 
& & + |fuf -v|_2 \\  
\\ & \leq & 4CK\epsilon + \sum_{1\leq i \leq p,  
1\leq j \leq k_i} (|y_{jj}^{(i)}uy_{jj}^{(i)} - 
y_{j1}^{(i)}uy_{1j}^{(i)}|_2 
+|y_{j1}^{(i)}uy_{1j}^{(i)} - y_{j1}^{(i)}v_iy_{1j}^{(i)}|_2) \\ 
& \leq & 4CK\epsilon + \sum_{1\leq i \leq p, 1\leq j \leq 
k_i}(2C\epsilon + |y_{11}^{(i)}uy_{11}^{(i)} - v_i|_2) \\ & 
\leq & 8CK\epsilon. \end{eqnarray*}

\noindent Set $L=8CK,$ observe that $L$ is dependent only on $x,$ 
and that $\inf_{h\in H_{\pi}}|u-h|_2 \leq |u-z|_2 \leq L \epsilon$.
\end{proof} 

	Denote $L > 0$ to be the uniform Lipschitz constant of Lemma 5.4.  
There exists a polynomial $f$ in $n$ noncommuting variables satisfying 
$f(a_1,\ldots,a_n)=x$.  There exists an $L_1 > 0$ such that for any $k\in 
\mathbb N$ and $ \xi_1,\ldots,\xi_n,\eta_1,\ldots,\eta_n \in 
(M_{k}^{sa}(\mathbb C))_R$ 

\[|f(\xi_1,\ldots,\xi_n) - f(\eta_1,\ldots,\eta_n)|_2 \leq L_1 \cdot 
\max\{|\xi_i-\eta_i|_2:1\leq i \leq n\}.\]
Denote $P_{\epsilon}(S)$ and $U(x_1,\ldots,x_d)$ to have the same meaning 
as in Section 4.

\begin{lemma}If $\rho_1, \rho_2 > 0, m\in \mathbb N$, and $\gamma>0$, then 
there is an $N\in \mathbb N$  such that for each $k \geq N$ there 
exists an $(x_{1}^{(k)},\ldots,x_{n}^{(k)}) \in 
\Gamma_R(a_1,\ldots,a_n;m,k,\gamma)$ satisfying for any 
$\frac{\kappa}{\rho_1 L L_1} > \epsilon > 0$

\[k^{-2} \cdot \log(P_{\epsilon 
\rho_1}(U(x_{1}^{(k)},\ldots,x_{n}^{(k)}))) \geq 
(\bigtriangleup_{\varphi}(M) - \rho_2) \cdot \log 
\left (\frac{\kappa}{\rho_1 L L_1 \epsilon} \right ). \]
\end{lemma}
\begin{proof}By Lemma 5.2 there is an $N \in \mathbb N$ such that for 
each $k \geq N$ there exists a $*$-homomorphism $\sigma_k:M \rightarrow 
M_k(\mathbb C)$ satisfying:
\begin{itemize}
\item $\|tr_k \circ \sigma_k - \varphi \| \leq \frac{\gamma}{(R+1)^m}$.
\item For each $k$ setting $H_k = H_{\sigma_k}$ and $X_k = X_{\sigma_k}$ 
we have that $H_k$ is a tractable 
Lie subgroup of $U_k$ and $k^2(\bigtriangleup_{\varphi}(M) - 
\rho_2) \leq \dim(X_k)$.
\item For $\epsilon > 0$

\[ \left(\frac{\kappa}{\epsilon}\right)^{\dim X_k} \leq
P(X_k, \epsilon)\]
where $P(X_k,\epsilon)$ is the maximum number of points in an $\epsilon$ 
separated subset of $X_k.$
\end{itemize} 
For each $1 \leq i \leq n$ and $k\geq N$ define $x_{i}^{(k)}=\sigma_k(a_i) 
\in M_{k}^{sa}(\mathbb C)$.  By the first condition above, 
$(x_{1}^{(k)},\ldots,x_{n}^{(k)}) \in 
\Gamma_R(a_1,\ldots,a_n;m,k,\gamma)$.  
To see that the second condition is fulfilled suppose $\epsilon$ satisfies 
the hypothesis of the lemma and $k \geq N$.  $\kappa \geq \rho_1LL_1 
\epsilon$ so  
\[P(X_k, \rho_1LL_1\epsilon) \geq \left( \frac 
{\kappa}{\rho_1LL_1\epsilon} \right)^{\dim X_k} 
\geq \left( \frac {\kappa}{\rho_1LL_1\epsilon} \right)^ 
{k^2(\bigtriangleup_{\rho}(M) - \rho_2)}.\]

\noindent For any $u,v \in U_k$

\begin{eqnarray}
d_{\sigma_k}(\dot{u},\dot{v}) & \leq & 
L \cdot |u \sigma_k(x)u^* - 
v \sigma_k(x)v^*|_2 \nonumber \\ & = & 
L \cdot 
|f(u \sigma_k(a_1) u^{*},\ldots,u \sigma_k(a_n)u^{*}) 
- 
f(v \sigma_k(a_1) v^*,\ldots,v \sigma_k(a_n) v^{*})|_2 
\nonumber \\ & \leq & LL_1 \cdot \max \{ 
|u x_{i}^{(k)}u^{*}- vx_{i}^{(k)}v^{*}|_2: 
1\leq i \leq n \}. \nonumber 
\end{eqnarray}
It follows that $P_{\epsilon 
\rho_1}(U(x_{1}^{(k)},\ldots,x_{n}^{(k)})) \geq P(X_k, 
\rho_1LL_1 
\epsilon).$  Hence for $k \geq N$
\[k^{-2} \cdot \log(P_{\epsilon
\rho_1}(U(x_{1}^{(k)},\ldots,x_{n}^{(k)}))) \geq
(\bigtriangleup_{\varphi}(M) - \rho_2) \cdot \log
\left (\frac{\kappa}{\rho_1 L L_1 \epsilon} \right ). \]
\end{proof}

	The following corollary will not be used until the next section.  
Suppose $r > R.$
	  
\begin{corollary} If $\Omega >0$ and $ \frac {\kappa}{3 \Omega L L_1} >
\epsilon > 0$, then there exist $m \in \mathbb N$ and $\gamma >0$ such
that if $\langle (x_1^{(k)}, \ldots,x_n^{(k)}) \rangle_{k=1}^{\infty}$
is a sequence satisfying $(x_1^{(k)}, \ldots,x_n^{(k)}) \in
\Gamma_r(a_1, \ldots, a_n; m, k, \gamma)$ for sufficiently large $k$,
then

\[ k^{-2} \cdot \log (P_{\Omega \epsilon} (U(x_1^{(k)}, \ldots,
x_n^{(k)})))  \geq (\bigtriangleup_{\varphi}(M) - \epsilon) \cdot \log
\left ( \frac {\kappa} { 3 L L_1 \Omega \epsilon} \right)\]

\noindent for sufficiently large $k.$
\end{corollary}

\begin{proof} By Corollary 3.4 there exist $m \in \mathbb N$ and $\gamma >0$ such 
that for any $k \in \mathbb N$ and $(y_1,\ldots,y_n)$, $(z_1,\ldots,z_n)$ $\in 
\Gamma_r(a_1, \ldots, a_n; m, k, \gamma)$ there exists a $u \in U_k$ such 
that for $1 \leq j \leq n$ $|uy_ju^{*} - z_j)|_2 < \frac {\epsilon}{2}.$  By 
Lemma 5.5 there exists an $N \in \mathbb N$ such that for each $k \geq N$ there 
exists an $(z_1^{(k)}, \ldots, z_n^{(k)}) \in \Gamma_r(a_1, \ldots, a_n; m, k, \gamma)$
satisfying 
 
\[ k^{-2} \cdot \log(P_{3 \Omega \epsilon}(U(z_1^{(k)},\ldots,z_n^{(k)})))
\geq (\bigtriangleup_{\varphi(M)} - \epsilon) \cdot \left ( \frac
{\kappa}{3 L L_1 \Omega \epsilon} \right). \]
 
\noindent Now merely observe that for any such sequence $\langle
(x_1^{(k)}, \ldots,x_n^{(k)})\rangle _{k=1}^{\infty}$ satisfying the
hypothesis of the corollary with $m$ and $\gamma$ chosen above,
 $P_{\Omega \epsilon}(U(x_1^{(k)}, \ldots,x_n^{(k)})) \geq P_{3 \Omega
\epsilon}(U(z_1^{(k)},\ldots,z_n^{(k)}))$ for $k \geq N.$ \end{proof}

\begin{theorem}$\delta_0(a_1,\ldots,a_n) \geq 1 - \sum_{i=1}^{p} \frac
{\alpha_{i}^{2}}{k_{i}^{2}}$. \end{theorem}

\begin{proof}Suppose $\min \left \{ \frac{1}{2}, \frac {\kappa}{4\sqrt{n}
LL_1} \right \} > \epsilon > 0, m \in \mathbb N$, and $\gamma, r > 0$.  
Corollary 2.14 of [11] provides an $N\in \mathbb N$ such that if $k \geq
N$ and $\sigma$ is a Radon probability measure on $((M_{k}^{sa}(\mathbb
C))_{R+1})^{2n}$ invariant under the $U_k$-action
$(\xi_1,\ldots,\xi_n,\eta_1,\ldots,\eta_n) \mapsto
(\xi_1,\ldots,\xi_n,u\eta_1 u^{*},\ldots,u\eta_n u^{*})$ where $ u\in
U_k$, then $\sigma(\omega_k) > \frac {1}{2}$ where

\begin{eqnarray}\omega_k = \{ (\xi_1,\ldots,\xi_n,\eta_1,\ldots,\eta_n)
\in((M_{k}^{sa}(\mathbb 
C))_{R+1})^{2n}\hspace{-.05in} 
& : & \hspace{-.05in}\{\xi_1,\ldots,\xi_n\} \text { and } 
\{\eta_1,\ldots,\eta_n\} \nonumber \\ & & \text { are } 
\left( m,\frac {\gamma}{4^{m}} \right)  
\text{ - free} \}. \nonumber
\end{eqnarray}
Lemma 5.5 provides an $N_1 \in \mathbb N$ such that for each $k \geq N_1$ 
there exists an $(x_{1}^{(k)},\ldots,x_{n}^{(k)}) \in  
\Gamma_R (a_1,\ldots,a_n;m,k, \gamma / (8(R+2))^m)$ 
satisfying

\[k^{-2} \cdot \log(P_{4\epsilon
\sqrt{n}}(U(x_{1}^{(k)},\ldots,x_{n}^{(k)}))) \geq
(\bigtriangleup_{\varphi}(M) - r) \cdot \log
\left (\frac{\kappa}{4 \sqrt{n} L L_1 \epsilon} \right ). \]
For $k \geq N+N_1$ denote by $\nu_k$ the 
atomic probability measure on $((M_{k}^{sa}(\mathbb C))_{R+1})^n$ 
concentrated at $(x_{1}^{(k)},\ldots,x_{n}^{(k)})$ and denote by $m_k$ the 
Radon probability measure obtained by restricting vol to
$\Gamma_{2\epsilon}(\epsilon s_1,\ldots,\epsilon s_n;
m,k, \gamma / 8^m)$ and normalizing appropriately.  $\nu_k \times 
m_k$ is a Radon probability measure on $((M_{k}^{sa}(\mathbb 
C))_{R+1})^{2n}$ invariant 
under the $U_k$-action described above.  Write $F_k$ for the set of all $(y_1,\ldots,y_n) \in \Gamma_{2\epsilon}(\epsilon
s_1,\ldots,\epsilon s_n;m,k,\gamma /8^m)$ such that  $\{y_1,\ldots,y_n\}$ 
and $\{x_{1}^{(k)},\ldots,x_{n}^{(k)}\}$ are $\left (m, \frac {\gamma}{4^m}\right)$-free.

For $k \geq N+N_1 \hspace{0.1in} \frac {1}{2} < (\nu_k \times 
m_k)(\omega_k) = 
m_k(F_k)$.  Set $E_k=(x_{1}^{(k)},\ldots,x_{n}^{(k)}) + F_k$.  For $k \geq 
N+N_1$ and $u\in U_k, \text {vol}(uE_ku^*) = \text {vol}(F_k)$ (where 
$uE_ku^*$ is defined as in Section 4) and $uE_ku^*$ is contained in 
$\Gamma_{R+1}(a_1 + \epsilon s_1,\ldots, a_n + \epsilon s_n: \epsilon 
s_1,\ldots, \epsilon s_n; m,k,\gamma).$

	For each $k \geq N+N_1$ find a subset $\langle u_{k,s}
\rangle_{s\in S_k}$ of $U_k$ such that $|S_k| =
P_{4\epsilon\sqrt{n}}(U(x_{1}^{(k)},\ldots,x_{n}^{(k)}))$ and for any
$s, s^{\prime} \in S_k, s \neq s^{\prime}$,

\[\max\{|u_{k,s}x_{i}^{(k)}u_{k,s}^{*} - 
u_{k,s^{\prime}}x_{i}^{(k)}u_{k,s^{\prime}}^{*}|_2: 1 \leq i \leq n \} > 
4\epsilon \sqrt{n}.\]
$F_k \subset (M_{k}^{sa}(\mathbb C))^n$ is a $2\epsilon\sqrt {nk}$ bounded 
subset with respect to the $\| \cdot \|_2$-norm.  Hence for any 
$s,s^{\prime} \in S_k, s\neq s^{\prime} (u_{k,s}E_ku_{k,s}^{*}) \bigcap 
(u_{k,s^{\prime}}E_ku_{k,s^{\prime}}^{*}) = \emptyset$.  Consequently for 
$k \geq N+N_1 \hspace{0.1in} \text{vol}(\Gamma_{R+1}(a_1 + \epsilon 
s_1,\ldots,a_n + 
\epsilon s_n: \epsilon s_1,\ldots, \epsilon s_n; m,k,\gamma))$ dominates

\begin{eqnarray} \text{vol}(\bigsqcup_{s\in S_k} u_{k,s}E_ku_{k,s}^{*}) =
|S_k| \cdot \text{vol}(F_k) & = & |S_k| \cdot m_k(F_k) \cdot
\text{vol}(\Gamma_{2\epsilon}(\epsilon s_1,\ldots,\epsilon s_n; 
m, k, \gamma /8^m )) \nonumber \\ & > & \frac {1}{2} 
\cdot |S_k| \cdot \text{vol}(\Gamma_{2\epsilon}(\epsilon 
s_1,\ldots,\epsilon s_n; m,k,\gamma / 8^m)). \nonumber\end{eqnarray}

	By what preceded for $\min \left \{ \frac{1}{2}, 
\frac{\kappa}{4\sqrt{n}LL_1} \right\}>\epsilon>0, m\in 
\mathbb N$, and $\gamma, r >0 \hspace{.1in} \chi_{R+1}(a_1+ \epsilon 
s_1,\ldots,a_n + 
\epsilon s_n: \epsilon s_1,\ldots,\epsilon s_n; m, \gamma)$ dominates 

\begin{eqnarray} & & \limsup_{k \rightarrow \infty}k^{-2} \cdot \log 
\left ( \frac{1}{2} \cdot |S_k| \cdot 
\text{vol}(\Gamma_{2\epsilon}(\epsilon 
s_1,\ldots,\epsilon s_n; m,k, \gamma /8^m)) + \frac{n}{2} \cdot 
\log k \right ) \nonumber \\ & = & \limsup_{k \rightarrow \infty} 
\left [ k^{-2} \cdot \log \left (\text{vol}(\Gamma_{2\epsilon}(\epsilon 
s_1,\ldots,\epsilon s_n; m,k, \gamma / 8^m)) + \frac{n}{2} \cdot  
\log k \right) + k^{-2} \cdot \log(|S_k|) \right] \nonumber \\ & \geq & 
\chi (\epsilon s_1,\ldots,\epsilon s_n) + 
\log \left (\left( \frac {\kappa}{4\sqrt{n}LL_1\epsilon} \right 
)^{\bigtriangleup_{\varphi}(M) - r} \right)  
\nonumber \\ & = & \log(\epsilon^{n+r-\bigtriangleup_{\varphi}(M)}) + 
\log \left ((2\pi e)^{\frac{n}{2}} \left (\frac {\kappa}{4\sqrt{n}LL_1} 
\right)^{\bigtriangleup_{\varphi}(M) - r} \right ). \nonumber\end{eqnarray} 
Letting $r \rightarrow 0$ it follows that 
\[\chi_{R+1}(a_1+ \epsilon s_1,\ldots, a_n + \epsilon s_n; m, \gamma) \geq 
\log(\epsilon^{n-\bigtriangleup_{\varphi}(M)}) + \log \left ((2\pi e)^{\frac{n}{2}} \left 
(\frac {\kappa}{4\sqrt{n}LL_1}
\right)^{\bigtriangleup_{\varphi}(M)} \right ).\]
This inequality holding for all $\epsilon > 0$ sufficiently small, 
$m\in \mathbb N$, and $\gamma > 0$
\[ \chi(a_1 + \epsilon s_1,\ldots,a_n + \epsilon s_n: s_1,\ldots,s_n) \geq 
\log(\epsilon^{n-\bigtriangleup_{\varphi}(M)}) + \log \left ((2\pi e)^{\frac{n}{2}} \left
(\frac {\kappa}{4\sqrt{n}LL_1}
\right)^{\bigtriangleup_{\varphi}(M)} \right ).\]
Dividing by $|\log \epsilon|$, taking lim sup's as $\epsilon$ goes to
$0$, and adding $n$ to both sides above yields 

\[\delta_0(a_1,\ldots,a_n) \geq 1 - \sum_{i=1}^{p} \frac {\alpha_{i}^{2}}{k_{i}^{2}}.\]
\end{proof}
By Theorem 3.10 and Theorem 5.7 we have:
\begin{corollary}$\delta_0(a_1,\ldots,a_n) = 1
- \sum_{i=1}^{p} \frac {\alpha_{i}^{2}}{k_{i}^{2}}$.
\end{corollary}

\section{The General Lower Bound}

	In this section we find a lower bound for $\delta_0(a_1,\ldots,a_n).$ 
When $M$ is hyperfinite the lower bound will be sharp.  By decomposing M over
its center it follows that 

\[M \simeq M_0 \oplus \hspace{.05in} (\oplus_{i=1}^{s}
M_{k_i}(\mathbb C)) \oplus M_{\infty}, \hspace{.1in} \varphi \simeq 
\alpha_0 
\varphi_0 \oplus \hspace {.05in} (\oplus_{i=1}^{s} \alpha_i 
tr_{k_i}) \oplus 0\] 

\noindent where as in the introduction $s \in \mathbb N \bigcup \{0\}
\bigcup \{\infty\}, \alpha_i > 0$ for $1 \leq i \leq s$ $(i \in \mathbb
N),$ $M_0$ is a diffuse von Neumann algebra or $\{0\}$, $\varphi_0$ is a
faithful, tracial state on $M_0$ and $\alpha_0 >0$ if $M_0 \neq \{0\}$,
$\varphi_0 = 0$ and $\alpha_0 =0$ if $M_0 = \{0\}$, and $M_{\infty}$ is
a von Neumann algebra or $\{0\}$. We remark that $M,$ hyperfinite or
otherwise, always admits such a decomposition.  We will show that
$\delta_0(a_1,\ldots,a_n) \geq 1 - \sum_{i=1}^{s} \frac
{\alpha_{i}^{2}}{k_{i}^{2}}.$ Again, because $\varphi$ vanishes on
$M_{\infty}$ and our main claim concerns the calculation of a lower
bound for $\delta_0(a_1,\ldots,a_n)$ assume without loss of generality
that

\[ M = M_0 \oplus (\oplus_{i=1}^{s} M_{k_i}(\mathbb C)), \varphi =
\alpha_0\varphi_0 \oplus (\oplus_{i=1}^{s} \alpha_i tr_{k_i}). \]

We proceed first by finding a suitable set of elements
$\{a_1^{\prime},a_2^{\prime},a_{3}^{\prime}\}$ in the $*$-algebra
generated by the $a_i$ such that the packing number of unitary orbits of
\emph{certain} microstates of
$\Gamma_C(a_1^{\prime},a_2^{\prime},a_3^{\prime};m,k,\gamma)$ approximate
(from below and in a normalized sense) $1 - \sum_{i=1}^{s} \frac
{\alpha_{i}^{2}}{k_{i}^{2}}$.  These microstates can be obtained as
noncommuting polynomials of well approximating microstates for
$\{a_1,\ldots,a_n\}$ and hence the packing number of unitary orbits of
such microstates of
$\Gamma_C(a_1^{\prime},a_2^{\prime},a_3^{\prime};m,k,\gamma)$ will provide
a lower bound for the packing number of unitary orbits of the microstates
for $\{a_1,\ldots,a_n\}.$ One can then use asymptotic freeness results to
transform these metric entropy quantities into free entropy dimension
quantities as in Theorem 4.5 and Theorem 5.7.

	Throughout this section write $A$ for the $*$-algebra generated by
$\{a_1,\ldots,a_n\}$.  We also maintain the notation introduced at the
beginning of Section 4.  In subsections 6.1 and 6.2 we assume that $M_0
\neq \{0\}$ (which implies $0 < \alpha_0$) and $\alpha_0 < 1.$

\subsection{Construction of
$\{a_{1}^{\prime},a_{2}^{\prime},a_{3}^{\prime}\}$ when $M_0 \neq \{0\}$
and $\alpha_0 <1$}

	Fix $l \in \mathbb N$ with $l \leq s$ and $\varepsilon >0$.  
Define $M_1=\oplus_{j=1}^{l} M_{k_j}(\mathbb C)$ and
$M_2=\oplus_{l<j \leq s} M_{k_j}(\mathbb C)$.  $M_1$ is a finite
dimensional C*-algebra and by Lemma 5.1 has two self-adjoint generators
$b_1$ and $b_2$.  $A$ is strongly dense in $M$ so $Ae$ is strongly dense
in $Me$ where $e = 0 \oplus (\oplus_{j=1}^{l} I_j) \oplus 0 \in
Z(M)$ and the $I_j$ are as in Section 3.  Thus $Ae = Me$.  Consequently
there exist $a_{1}^{\prime}, a_{2}^{\prime} \in A$ such that
\[a_{i}^{\prime}=f_i \oplus b_i \oplus \xi_i \in M_0 \oplus M_1
\oplus M_2 = M.\] $M_0$ being diffuse, there exists an $f \in M_0$ such
that $\delta_0(f)  =1$ and $sp(f)=[1,2]$ (here $\delta_0(f)$ is calculated
with respect to $\varphi_0$).  $Ae_0$ is strongly dense in $Me_0$ where
$e_0 = I_0 \oplus 0 \in M_0 \oplus (\oplus_{1 \leq j \leq s}
M_{k_j}(\mathbb C).$ By Kaplansky's Density Theorem, Proposition 6.14 of
[8], and Corollary 6.7 of [9] there exists an $0 \leq a \in A$ satisfying:
\begin{itemize} \item $a=g \oplus b \oplus \xi \in M_0 \oplus M_1
\oplus M_2 = M$. \item $sp(g) \subset [0,2],\hspace{.08in} \delta_0(g)
> 1 - \varepsilon$, and $ \varphi_0(\chi_{[0, 1/2]}) < \varepsilon$ where
for any Borel subset $S \subset \mathbb R \hspace{.08in} \chi_S$ denotes
the spectral projection of $g$ associated to the set $S$. \item $\|b\| <
\frac{1}{6}$. \end{itemize}

Since $b \geq 0$ and $M_1$ is finite dimensional 
$sp(b)=\{\beta_1,\ldots,\beta_d\} \subset \mathbb R$ where $0 \leq \beta_1 
< \ldots < \beta_d \leq \frac {1}{6}.$  Define $h: [0,2] \rightarrow \mathbb R$ by 
$h(t)=(t-\beta_1) \cdots (t-\beta_d)$.  Observe that:
\begin{itemize}
\item $h(b) = 0$.
\item $h^{-1}(h([0, 2\beta_d])) \subset [0,3\beta_d] \subset [0, 1/2]$.
\item $h^{-1}(h((2\beta_d,2])) \subset (\beta_d,2]$.
\end{itemize}
The third observation implies that if $\beta \in h(2\beta_d, 2]$ then 
$h^{-1}({\beta})$ consists of exactly one point in $(\beta_d,2]$ since $h$ 
is strictly 
increasing on $(\beta_d,2]$.  Noting that for all but countably many 
$\beta \hspace{.1in} \chi_{h^{-1}\{\beta\}} = 0$

\begin{eqnarray}
\sum_{\beta \in sp (h(b))} \varphi_0(\chi_{h^{-1}\{\beta\}})^2 & = & 
\sum_{\beta \in h(sp(b))} \varphi_0(\chi_{h^{-1}\{\beta\}})^2 
\nonumber \\ & \leq & \sum_{\beta \in h([0, 2\beta_d])} 
\varphi_0(\chi_{h^{-1}\{\beta\}})^2  + \sum_{\beta \in h((2\beta_d, 2])} 
\varphi_0(\chi_{h^{-1}\{\beta\}})^2 \nonumber \\ & \leq & 
\varphi_0(\sum_{\beta \in h([0,2\beta_d])} \chi_{h^{-1}\{\beta\}} ) + 
\sum_{\beta \in [0,2]} \varphi_0(\chi_{\{\beta\}})^2 \nonumber \\  & 
\leq & 
\varphi_0(\chi_{[0,3\beta_d]}) + 
\varepsilon 
\nonumber \\ & < &
2\varepsilon. \nonumber
\end{eqnarray}
Define $a_{3}^{\prime} = h(a) = h(g) \oplus 0 \oplus h(\xi)$.  We 
have just proven:
\begin{lemma} If $l\in \mathbb N, 1 \leq l\leq s,$ and $\varepsilon > 0$, 
then there exist $a_{1}^{\prime}, a_{2}^{\prime}, a_{3}^{\prime} \in A $ 
of the form $a_{i}^{\prime} = f_i \oplus b_i \oplus \xi_i \in M_0 
\oplus M_1 \oplus M_2 = M$ satisfying:
\begin{itemize}
\item \{$b_1$, $b_2$\} generates $M_1$ 
and $b_3=0$.
\item $\delta_0(f_3) > 1- \varepsilon$.
\end{itemize}
\end{lemma} 

\subsection{Lower Bounds Estimates for $\delta_0(a_1,\ldots,a_n)$ when $M_0 \neq \{0\}$ and $\alpha_0 <1$}

	Fix $l\in \mathbb N$, $ l \leq s$, and $\varepsilon >0.$ Find
elements $a_{1}^{\prime}, a_{2}^{\prime}, a_{3}^{\prime} \in A$ with the
properties listed in Lemma 6.1 and denote $M_1, M_2$ to have the same
meaning as in the preceding subsection.  Suppose $C \geq
\max\{\|a_{i}^{\prime}\|\}_{i=1,2,3}+1$.  Denote by $\varphi_1$ the
tracial state on $M_1$ obtained by restricting $\varphi$ to $0 \oplus
M_1 \oplus 0$ and normalizing appropriately.  Similarly denote by
$\varphi_2$ the positive trace (possibly 0) on $M_2$ obtained by
restricting $\varphi$ to $0 \oplus 0 \oplus M_2$ and normalizing
appropriately.  Define $\beta_0=\alpha_0, \beta_1= \sum_{i=1}^{l}
\alpha_i$, $\beta = \min \{\beta_0, \beta_1\},$ and $\beta_2 = 1- \beta_0
- \beta_1.$ Recall the constant $\kappa, L,$ and $L_1$ of the previous
section with respect to $M_1$, $\varphi_1$, and $x = b_1 + i b_2$.  
Finally define $e_0$ and $e_1$ to be the projections $I_0 \oplus 0
\oplus 0$, $0 \oplus (\oplus_{j=1}^{l} I_j) \oplus 0 \in M$,
respectively, and $e_2 = I -e_0 - e_1.$

\begin{lemma} If $D >0$ and $\min \left \{\frac {\beta \kappa}{27DLL_1},
1 \right \} > \epsilon >0$, then there exist $m_{\epsilon} \in \mathbb
N$ and $\gamma_{\epsilon} >0$ (dependent on $\epsilon$) such if $\langle
(z_1^{(k)},z_2^{(k)},z_3^{(k)}) \rangle_{k=1}^{\infty}$ is a sequence
satisfying $(z_1^{(k)},z_2^{(k)},z_3^{(k)}) \in (M_{k}^{sa}(\mathbb
C))^3$ for all $k$ and $(z_1^{(k)},z_2^{(k)},z_3^{(k)}) \in
\Gamma_C(a_1^{\prime},a_2^{\prime},a_3^{\prime} : e_0,e_1,e_2;
m_{\epsilon},k,\gamma_{\epsilon})$ for sufficiently large $k$, then
$\limsup_{k \rightarrow \infty} k^{-2} \cdot \log (P_{D \epsilon}
(U(z_1^{(k)},z_2^{(k)},z_3^{(k)})))$ dominates

\[
(\beta_0 + \beta_1 - \epsilon)^2 \cdot
(\chi((f_3 \oplus 0 ) + \epsilon s: s) + |\log \epsilon| - K) +
(\beta_1 - \epsilon)^2(\bigtriangleup_{\varphi_1}(M_1) - \epsilon)\cdot
\log \left(\frac {\beta \kappa}{27 \epsilon D LL_1}\right)
\]
where $K=\log((2 + 9 \beta^{-1} D) \sqrt{2\pi e})$ and $s$ is a 
semicircular element free with respect to $M_0 \oplus M_1.$

\end{lemma}

\begin{proof}

	Suppose $\min \left \{\frac {\beta \kappa}{27 DLL_1}, 1 \right \}
> \epsilon >0.$ By Lemma 4.3 there exist an $m_1 \in \mathbb N$ and
$\gamma_1$ (dependent on $\epsilon$) such that if
$\langle x^{(k)} \rangle_{k=1}^{\infty}$ is a sequence satisfying 
$x^{(k)} \in
M_{k}^{sa}(\mathbb C)$ for all $k\in \mathbb N$ and $x^{(k)} \in
\Gamma_C(f_3 \oplus 0; m_1,k,\gamma_1)$ for sufficiently large $k$,
then 

\[ \limsup_{k\rightarrow \infty} [ k^{-2} \cdot \log P_{9 \beta^{-1}
D \epsilon}(U(x^{(k)}))] \geq \chi_C((f_3 \oplus 0) + \epsilon s :s) +
|\log \epsilon| - K\] 

\noindent where $K = \log ((2+ 9 \beta^{-1} D)\sqrt{2\pi e}), f_3
\oplus 0 \in M_0 \oplus M_1$, $M_0 \oplus M_1$ is endowed with
the tracial state $(\beta_0 + \beta_1)^{-1} (\beta_0 \varphi_0 \oplus
\beta_1 \varphi_1)$, and $s$ is a semicircular element free with respect
to $M_0 \oplus M_1$.
	
	By Corollary 5.6 there exist an $m_2 \in \mathbb N$ and $\gamma_2 >0$
such that if  $\langle (y_1^{(k)}, y_2^{(k)})\rangle_{k=1}^{\infty}$ is 
a sequence 
satisfying $(y_1^{(k)}, y_2^{(k)}) \in (M_k^{sa}(\mathbb C))^2$ for all $k$ 
and $(y_1^{(k)}, y_2^{(k)}) \in \Gamma_C(b_1,b_2;m_2,k, \gamma_2)$ for 
sufficiently large $k$, then
	
\[ k^{-2} \cdot \log P_{9 \beta^{-1} D \epsilon} (U(y_1^{(k)}, y_2^{(k)}))
\geq \log \left( \left( \frac {\beta \kappa}{27 \epsilon D L L_1} \right
)^{\bigtriangleup_{\varphi_1}(M_1) - \epsilon} \right) \]

\noindent for sufficiently large $k$.

	By standard approximations for any $\varepsilon >0$ there exist $m \in 
\mathbb N$ and $\gamma >0$ such that for any $k \in \mathbb N$ if $(z_1,z_2,z_3) \in 
\Gamma_C(a_1^{\prime},a_2^{\prime},a_3^{\prime}:e_0,e_1,e_2; m, 
k, \gamma ),$ then the following conditions are satisfied:

\begin{itemize}

\item There exist mutually orthogonal projections $q_0,q_1,q_2 \in 
M_k(\mathbb C)$ with $n_0+n_1+n_2 = k$ where $n_i$ denotes the 
dimension of the range of $q_i$ and for each $i$ $|tr_k(q_i) - 
\beta_i| < \varepsilon.$

\item Canonically identifying $q_0M_k(\mathbb C)q_0 + q_1M_k(\mathbb C)q_1 
+ q_2M_k(\mathbb C)q_2$ with  $\oplus_{i=0}^2 M_{n_i}(\mathbb C)$
$|z_i - h_i|_2 < \varepsilon $ where $h_i = x_i \oplus y_i \oplus
\eta_i \in \oplus_{i=0}^2 M_{n_i}(\mathbb C) \subset M_k(\mathbb C)$ and $y_3 = 0.$   

\item $ x_3 \oplus 0  \in \Gamma_C(f_3 \oplus 0; m_1, n_0+n_1, 
\gamma_1).$

\item $(y_1,y_2) \in \Gamma_C(b_1,b_2;m_2, n_1,\gamma_2)$ where $b_1,
b_2 \in (M, \varphi_1).$ \end{itemize}

\noindent Secondly given $\varepsilon >0$ there exists an $N \in \mathbb
N$ such that for each $k > N$ there exists a corresponding $k < \tau(k)
\in \mathbb N$ satisfying $| \frac {k}{\tau(k)} - (\beta_0 + \beta_1)| <
\varepsilon.$ Combining these two remarks there exist $m_{\epsilon} \in
\mathbb N$ and $\gamma_{\epsilon} >0$ such that for any $k \in \mathbb N$
sufficiently large there exists a $k < \tau(k) \in \mathbb N$ such that if
$(z_1,z_2,z_3) \in
\Gamma_C(a_1^{\prime},a_2^{\prime},a_3^{\prime}:e_0,e_1,e_2; m_{\epsilon},
\tau(k), \gamma_{\epsilon})$, then the following conditions are satisfied:

\begin{itemize}

\item There exist mutually orthogonal projections $q_0,q_1,q_2 \in 
M_{\tau(k)}(\mathbb C)$ with $n_0+n_1+n_2 = \tau(k)$ where $n_i$ denotes the 
dimension of the range of $q_i$ and for each $i$ $|tr_{\tau(k)}(q_i) - 
\beta_i| < \frac {1}{2} \cdot \min \{ \epsilon, \beta \}$.     

\item $n_0 +n_1 = k$
 
\item Canonically identifying $q_0M_{\tau(k)}(\mathbb C)q_0 +
q_1M_{\tau(k)}(\mathbb C)q_1 + q_2M_{\tau(k)}(\mathbb C)q_2$ with
$\oplus_{i=0}^2 M_{n_i}(\mathbb C)$ $|z_i - h_i|_2 < \frac
{\epsilon}{2}$ where $h_i = x_i \oplus y_i \oplus \eta_i \in
\oplus_{i=0}^2 M_{n_i}(\mathbb C) \subset M_{\tau(k)}(\mathbb C)$ and
$y_3 = 0$.

\item $ x_3 \oplus 0  \in \Gamma_C(f_3 \oplus 0; m_1, k, 
\gamma_1)$.

\item $(y_1,y_2) \in \Gamma_C(b_1,b_2;m_2, n_1,\gamma_2)$ where $b_1, b_2
\in (M, \varphi_1).$ \end{itemize}

	Now fix $k$ and $(z_1,z_2,z_3)$ satisfying the aforementioned
conditions so that the five properties listed just above hold.  Let
$n_i, h_i, x_i, y_i,$ and $\eta_i$ correspond to the fixed $k$ and
$(z_1,z_2,z_3).$ Find a set of unitaries $\langle u_s\rangle_{s \in
S_k}$ of $U_k$ suchthat $|S| = P_{9 \beta^{-1} D \epsilon}(U(x_3
\oplus 0))$ and for any $s, s^{\prime} \in S$ and $s \neq
s^{\prime}$,

\[ |u_s(x_3 \oplus 0)u_{s}^{*} - u_{s^{\prime}}(x_3 \oplus
0)u_{s^{\prime}}^{*}|_2 > 9 \beta^{-1} D \epsilon.\]

\noindent Find a set of unitaries $\langle v_g \rangle_{g \in G}$ of
$U_{n_1}$ such that $|G| = P_{9 \beta^{-1} D \epsilon}(U(y_1,y_2))$ and
for any $g, g^{\prime} \in G$ and $g \neq g^{\prime}$ 

\[ \max\{|v_gy_iv_{g}^{*} - v_{g^{\prime}}y_iv_{g_{\prime}}^{*}|_2 :  
i=1,2\} > 9 \beta^{-1} D \epsilon.\] 

\noindent For $(s,g) \in S_k \times G_k$ define $w_{s,g} \in U_{\tau(k)}$
by $w_{s,g}=(u_s(I_{n_0} \oplus v_g)) \oplus I_{n_2} \in
U_{\tau(k)}$.

	I claim that the family $\langle
(w_{s,g}h_1w_{s,g}^{*},w_{s,g}h_2w_{s,g}^{*},w_{s,g}h_3w_{s,g}^{*})\rangle_{s,g}
\in S\times G$ is a $3 D \epsilon $-separated set (with respect to the
$\rho$ metric defined in section 4).  Suppose
$(s,g),(s^{\prime},g^{\prime}) \in S \times G$ and $(s,g) \neq
(s^{\prime},g^{\prime})$.  Then either $s \neq s^{\prime}$ or $g \neq
g^{\prime}$.  In the former

\begin{eqnarray}|w_{s,g}h_3w_{s,g}^*-w_{s^{\prime},g^{\prime}}h_3w_{s^{\prime},g^{\prime}}^{*}|_2
& = & |w_{s,g}(x_3 \oplus 0 \oplus \eta_3)w_{s,g}^*-
w_{s^{\prime},g^{\prime}}(x_3 \oplus 0 \oplus
\eta_3)w_{s^{\prime},g^{\prime}}^{*}|_2 \nonumber \\ & = & |u_{s}(x_3
\oplus 0)u_{s}^{*} \oplus 0 - u_{s^{\prime}}(x_3 \oplus
0)u_{s^{\prime}}^{*} \oplus 0|_2 \nonumber \\ & > & \sqrt{ \frac
{\beta_1 + \beta_2}{2} } \cdot 9 \beta^{-1} D \epsilon \nonumber \\ & >
& 3 D \epsilon \nonumber.  \end{eqnarray}

\noindent Suppose $g \neq g^{\prime}$.  We can assume $s = s^{\prime}$.  
For $i=1,2$ 

\begin{eqnarray} |w_{s,g}h_iw_{s,g}^{*} -
w_{s,g^{\prime}}h_iw_{s,g^{\prime}}^{*}|_2 & = & |w_{s,g}(x_i \oplus
y_i \oplus \eta_i) w_{s,g}^{*} - w_{s, g^{\prime}}(x_i \oplus y_i
\oplus \eta_i)w_{s, g^{\prime}}^{*}|_2 \nonumber \\ & = & |u_s(x_i
\oplus v_g y_i v_{g}^{*})u_{s}^{*} \oplus \eta_3 - u_s(x_i \oplus
v_{g^{\prime}}y_iv_{g^{\prime}}^{*})u_{s}^{*} \oplus \eta_3|_2
\nonumber \\ & \geq & \sqrt{ \frac {\beta_1}{2}} \cdot
|v_gy_iv_{g}^{*}-v_{g^{\prime}}y_iv_{g^{\prime}}^{*}|_2 \nonumber
\end{eqnarray} 

\noindent so that 

\begin{eqnarray}\max \{|w_{s,g}h_iw_{s,g}^{*} -
w_{s,g^{\prime}}h_iw_{s,g^{\prime}}^{*}|_2 : i = 1,2 \} & \geq & \sqrt
{\frac {\beta_1}{2}} \cdot \max \{|v_gy_iv_{g}^{*} -
v_{g^{\prime}}y_iv_{g^{\prime}}^{*}|_2 : i=1,2\} \nonumber \\ & > & \sqrt
{\frac{\beta_1}{2}} \cdot 9 \beta^{-1} D \epsilon \nonumber\\ & > & 3 D
\epsilon \nonumber. \end{eqnarray}

\noindent By the inequalities above $P_{D \epsilon} 
(U(z_1,z_2,z_3)) \geq P_{3 D \epsilon} (U(h_1,h_2,h_3)) \geq |S \times 
G|.$  
	
	Now suppose $\langle (z_1^{(k)}, z_2^{(k)},z_3^{(k)})\rangle
_{k=1}^{\infty}$ is a sequence satisfying the hypothesis of the lemma
with $m_{\epsilon}$ and $\gamma_{\epsilon}$ as chosen on the previous
page.  For $k$ sufficiently large \[
(z_1^{(\tau(k))},z_2^{(\tau(k))},z_3^{(\tau(k))}) \in
\Gamma_C(a_1^{\prime}, a_2^{\prime},a_3^{\prime} : e_0,e_1,e_2;
m_{\epsilon}, \tau(k), \gamma_{\epsilon}).\] Thus
$(z_1^{(\tau(k))},z_2^{(\tau(k))},z_3^{(\tau(k))})$ satisfies the five
conditions previously stated.  For each $k$ sufficiently large consider
the corresponding $\tau(k) \in \mathbb N$ and denote $n_i(k)$,
$x^{(n_i(k))}$, $y_i^{(n_i(k))}$, $S_k$, and $G_k$ to be the $n_i$,
$x_i$, $y_i$, $S$, and $G$, respectively, associated to
$(z_1^{(\tau(k))},z_2^{(\tau(k))},z_3^{(\tau(k))}).$ \begin{eqnarray}
\limsup_{k \rightarrow \infty} \tau(k)^{-2} \cdot \log (P_{D \epsilon}
(U(z_1^{(\tau(k))},z_2^{(\tau(k))},z_3^{(\tau(k))})))  & \geq &
\limsup_{k \rightarrow \infty} \tau(k)^{-2} \cdot (\log |S_k| + \log
|G_k|) \nonumber \\ & \geq & \limsup_{k \rightarrow \infty}
(\tau(k)^{-2} \cdot \log |S_k|) + \nonumber \\ & & \liminf_{k
\rightarrow \infty} (\tau(k)^{-2} \cdot \log |G_k|) \nonumber.  
\end{eqnarray}

\noindent Set $\mathcal L = 9 \beta^{-1} D \epsilon.$ $\limsup_{k
\rightarrow \infty} (\tau(k)^{-2} \cdot \log |S_k|)$ dominates

\begin{eqnarray*} & & (\beta_0 + \beta_1 - \epsilon)^2 \cdot \limsup_{k
\rightarrow \infty} \frac {1}{(n_0(k)+n_1(k))^2} \cdot \log (P_{\mathcal
L}(U(x_3^{(n_0(k))} \oplus 0))) \nonumber \\ & = & (\beta_0 + \beta_1 -
\epsilon)^2 \cdot \limsup_{k \rightarrow \infty} \frac {1}{k^2} \cdot \log
(P_{\mathcal L}(U(x_3^{(n_0(k))} \oplus 0)))  \nonumber \\ & \geq &
(\beta_0 + \beta_1 - \epsilon)^2 \cdot \chi_C((f_3 \oplus 0) + \epsilon
s :s) + |\log \epsilon| - K \nonumber \end{eqnarray*} \noindent and since
$n_1(k) \rightarrow \infty$ as $k \rightarrow \infty$

\begin{eqnarray} \liminf_{k \rightarrow \infty} (\tau(k)^{-2} \cdot \log
|G_k|) & \geq & (\beta_1 - \epsilon)^2 \cdot \liminf_{k \rightarrow
\infty} \hspace{.1in} [n_1(k)^{-2} \cdot \log (P_{\mathcal L}
(U(y_1^{(n_1(k))},y_2^{(n_1(k))})))] \nonumber \\ & \geq & (\beta_1 -
\epsilon)^2 \cdot \log \left( \left( \frac {\beta \kappa}{27 \epsilon D L
L_1} \right )^{\bigtriangleup_{\varphi_1}(M_1)
 - \epsilon} \right). \nonumber
\end{eqnarray} 

\noindent The desired result follows.
\end{proof}

	We now recycle a familiar argument.  There exist polynomials $p_1,
p_2, p_3$ in $n$ noncommuting variables such that for $1 \leq j \leq 3$
$p_j(a_1, \ldots, a_n) = a_j^{\prime}$.  Without loss of generality we can
assume the $p_j$ take n-tuples of self-adjoint elements to self-adjoint
elements.  Find a constant $1 < \mathcal C$ such that that for any $(x_1,
\ldots, x_n) \in ((M_k^{sa}(\mathbb C))_{R+1})^n$ $\|p_j(x_1, \ldots, x_n)
\| < \mathcal C.$ Also there exists a $D_0 >0$ such that if $k \in \mathbb
N$ and $x_1, \ldots, x_n, y_1, \ldots, y_n \in (M_k^{sa}(\mathbb
C))_{R+1}$, then for all $j$

\[ |p_j(x_1, \ldots, x_n) - p_j(y_1, \ldots, y_n)|_2 \leq D_0 \cdot \max \{|x_i - y_i|_2 : 1 \leq i \leq n\}. \]

\begin{lemma} For $m \in \mathbb N,$ $\gamma >0,$ and small enough
$\epsilon > 0$ there is a sequence $\langle (x_1^{(k)},
\ldots,x_n^{(k)}) \rangle_{k=1}^{\infty}$ such that $(x_1^{(k)},
\ldots,x_n^{(k)}) \in (M_k^{sa}(\mathbb C))^n$ for all $k,$
$(x_1^{(k)},\ldots, x_n^{(k)}) \in \Gamma_{R+1}(a_1, \ldots, a_n ; m, k,
\gamma)$ for sufficiently large $k,$ and $\limsup_{k \rightarrow \infty}
k^{-2} \cdot \log (P_{4 \epsilon \sqrt {n} }(U(x_1^{(k)},
\ldots,x_n^{(k)})))$ dominates

\[(\beta_0 + \beta_1 - \epsilon)^2 \cdot
(\chi((f_3 \oplus 0 ) + \epsilon s: s) + |\log \epsilon| - K) +
(\beta_1 - \epsilon)^2(\bigtriangleup_{\varphi_1}(M_1) - \epsilon)\cdot
\log \left(\frac {\beta \kappa}{108 \sqrt{n} \epsilon D_0 LL_1}\right)
\]

\noindent where $K=\log((2 + 36 \beta^{-1} \sqrt{n} D_0) \sqrt{2\pi
e})$. \end{lemma} \begin{proof} Suppose $ \min \left \{ \frac {\beta
\kappa}{32 D_0 \sqrt{n} LL_1 }, 1 \right \} > \epsilon >1.$ By Lemma 6.2
there exist $m_{\epsilon} \in \mathbb N$ and $\gamma_{\epsilon} >0$ such
that if $\langle (z_1^{(k)},z_2^{(k)},z_3^{(k)}) \rangle_{k=1}^{\infty}$
is a sequence satisfying $(z_1^{(k)},z_2^{(k)},z_3^{(k)}) \in
(M_k^{sa}(\mathbb C))^3$ for all $k$ and
$(z_1^{(k)},z_2^{(k)},z_3^{(k)})  \in
\Gamma_C(a_1^{\prime},a_2^{\prime},a_3^{\prime}: e_1,e_2,e_3;
m_{\epsilon}, k,\gamma_{\epsilon})$ for sufficiently large $k$, then
$\limsup_{k\rightarrow \infty} k^{-2} \cdot \log (P_{4 \epsilon
\sqrt{n}D_0} (U(z_1^{(k)},z_2^{(k)},z_3^{(k)})))$ dominates

\[(\beta_0 + \beta_1 - \epsilon)^2 \cdot
(\chi((f_3 \oplus 0 ) + \epsilon s: s) + |\log \epsilon| - K) +
(\beta_1 - \epsilon)^2(\bigtriangleup_{\varphi_1}(M_1) - \epsilon)\cdot
\log \left(\frac {\beta \kappa}{108 \sqrt{n} \epsilon D_0 
LL_1}\right).\] 

  Because $\{a_1,\ldots,a_n\}$ has finite dimensional approximants, 
there exists an $N \in \mathbb N$ such that for each $k \geq N$ there is 
an $(x_1^{(k)},\ldots,x_n^{(k)}) \in ((M_k^{sa}(\mathbb C))_{R+1})^n$ 
satisfying 

\[(x_1^{(k)},\ldots,x_n^{(k)}) \in 
\Gamma_{R+1}(a_1,\ldots,a_n;m,k,\gamma)\]

\noindent and 

\[(y_1^{(k)}, y_2^{(k)},y_3^{(k)}) \in 
\Gamma_{\mathcal C}(a_1^{\prime},a_2^{\prime},a_3^{\prime}:e_1,e_2,e_3;
m_{\epsilon},k,\gamma_{\epsilon})\]

\noindent where for any $1 \leq j \leq 3$ $p_j(x_1^{(k)},\ldots,x_n^{(k)}) 
= y_j^{(k)}.$  Note that we can use the cutoff constant $\mathcal C$ 
by the argument of Lemma 4.1.

	If $u,v \in U_k$, then $D_0 \cdot \max 
\{|ux_i^{(k)}u^{*} -
vx_i^{(k)}v^{*}|_2 :  1 \leq i \leq n \}$ is greater than 

\[ |p_j(ux_1^{(k)}, 
\ldots,ux_n^{(k)}u^*) -
p_j(vx_1^{(k)}v^*, \ldots,vx_n^{(k)}v^*)|_2 = 
|uy_j^{(k)}u^* - vy_j^{(k)}v^*|_2. \]

\noindent Hence, 
\[\limsup_{k \rightarrow \infty}k^{-2} \cdot \log (P_{4 \epsilon \sqrt{n}}(U(x_1^{(k)},\ldots,x_n^{(k)}))) \geq \limsup_{k \rightarrow \infty} k^{-2} \cdot \log (P_{4\epsilon \sqrt{n} D_0}(U(y_1^{(k)},y_2^{(k)},y_3^{(k)}))).\]

\noindent In turn the dominated term is greater than or equal to

\[
(\beta_0 + \beta_1 - \epsilon)^2 \cdot
(\chi((f_3 \oplus 0 ) + \epsilon s: s) + |\log \epsilon| - K) +
(\beta_1 - \epsilon)^2(\bigtriangleup_{\varphi_1}(M_1) - \epsilon)\cdot
\log \left(\frac {\beta \kappa}{108 \sqrt{n} \epsilon D_0 LL_1}\right).
\]
\end{proof}

	 Lemma 6.3 more or less gives the lower bound.  We simply 
use the same sort of arguments which allowed us passage from
Lemma 4.4 to Theorem 4.5 and Lemma 5.5 to Theorem 5.7.  Namely
Voiculescu's approximate freeness results produce an $\epsilon$ ball, most
of whose elements are semicircular microstates trying very hard to be free
with respect to a matricial microstate for $\{a_1,\ldots,a_n
\}$.  Adding the $\epsilon$ ball to the single microstate creates
microstates for $\{a_1 + \epsilon s_1, \ldots ,a_n + \epsilon s_n\}$.  
The $\epsilon$ packing number of the unitary orbit of the microstate
yields the same number of disjoint, conjugate balls which are microstates
of $\{a_1 + \epsilon s_1, \ldots, a_n + \epsilon s_n\}$.  Applying $\log$
to this packing number, multiplying by $k^{-2}$ and taking a $\limsup$ as
$k \rightarrow \infty,$ dividing by $| \log \epsilon |$ and taking a
$\limsup$ as $\epsilon \rightarrow 0$ yields a lower bound for the
modified free entropy dimension of the n-tuple $\{a_1, \ldots,a_n \}.$

	By the discussion above $\delta_0(a_1,\ldots,a_n)$ should 
dominate the number obtained by taking the majorized quantity of Lemma 
6.3, dividing by $|\log \epsilon|$, and taking a $\limsup$ as $\epsilon 
\rightarrow 0$.  The resultant quantity of these successive operations is:

\begin{eqnarray*} (\beta_0 + \beta_1)^2( \delta_0(f_3 \oplus 0)) +
\beta_{1}^{2}\cdot
\bigtriangleup_{\varphi_1}(M_1) & > & (\beta_0 + \beta_1)^2 \cdot [1 -
(\beta_0 + \beta_1)^{-2}(\varepsilon \beta_0^2 + \beta_1^2)] \\ & & + 
\beta_{1}^{2}\cdot   
\bigtriangleup_{\varphi_1}(M_1) \nonumber \\
& > & (\beta_0 + \beta_1)^2 - ( \varepsilon + \beta_{1}^{2}) +
\beta_{1}^{2}\cdot \left (1 - \sum_{i=1}^{l}
\frac{\beta_{1}^{-2}\alpha_{i}^{2}}{k_{i}^{2}} \right)
\nonumber \\ & = & (\beta_{0} + \beta_{1})^2 - \varepsilon -
\sum_{i=1}^{l}\frac{\alpha_{i}^{2}}{k_{i}^{2}}. \nonumber
\end{eqnarray*}

	Omitting the details of a familiar analysis we conclude:

\begin{theorem}If $M_0 \neq \{0\}$ and $\alpha_0 <1,$ then
$\delta_0(a_1,\ldots,a_n)  \geq (\beta_{0} + \beta_{1})^2 - \varepsilon
-\sum_{i=1}^{l}\frac{\alpha_{i}^{2}}{k_{i}^{2}}.$ \end{theorem}

	Letting $l \rightarrow s$ and $\varepsilon \rightarrow 0$ we
arrive at:

\begin{corollary} If $M_0 \neq \{0\}$ and $\alpha_0 <1,$ then
$\delta_0(a_1,\ldots,a_n) \geq 1 - \sum_{i=1}^{s}
\frac{\alpha_{i}^{2}}{k_{i}^{2}}$.  \end{corollary}

\subsection{General Lower Bound Estimates for $\delta_0(a_1,\ldots,a_n)$}

	Now we deal with the situation where $M_0 = \{0\}$ or $\alpha_0
=1.$ Suppose first that $M_0 = \{0\}$.  The lower estimates come even
easier for in this case $M$ is merely a direct sum of matrix algebras.  
As in section 6.1 given $l \in \mathbb N$ with $l \leq s$ and using the
same notation, we can construct two self-adjoint elements $a_i^{\prime}
= b_i \oplus \xi_i \in M_1 \oplus M_2$ for $i =1,2$ such that they
lie in the $*$-algebra generated by $\{a_1,\ldots,a_n\}$ and $b_1, b_2$
generate $M_1$.  We have a similar packing number estimate in this case:

\begin{lemma}For small enough $\epsilon, \gamma>0$ and any $m\in \mathbb
N$ there is a sequence $\langle
(z_{1}^{(\bigtriangleup_k)},z_{2}^{(\bigtriangleup_k)})\rangle
_{k=1}^{\infty}$ such that $k < \bigtriangleup_k \in \mathbb N$ for all
$k$, $(z_{1}^{(\bigtriangleup_k)},z_{2}^{(\bigtriangleup_k)}) \in
M_{\bigtriangleup_k}^{sa}(\mathbb C)$ for all $k$,

\[(z_{1}^{(\bigtriangleup_k)},z_{2}^{(\bigtriangleup_k)}) 
\in 
\Gamma_C(a_{1}^{\prime},a_{2}^{\prime}; 
m, \bigtriangleup_k, \gamma)\]
for sufficiently large $k$, and 

\[\bigtriangleup_{k}^{-2}\cdot \log(P_{16  
\epsilon}(U(z_{1}^{(\bigtriangleup_k)},z_{2}^{(\bigtriangleup_k)}))) \geq
(\beta_1 - \epsilon)^2(\bigtriangleup_{\varphi_1}(M_1) - \epsilon)\cdot
\log \left(\frac {\kappa}{16 \beta^{-1} \epsilon LL_1}\right)
\]
for sufficiently large $k$.
\end{lemma}

The proof mimics that of Lemma 6.2 but it's even easier.  Lemma 6.6 requires 
that we find good enough microstates for $a_1^{\prime},a_2^{\prime}$ and this
is easy to do (in light of Lemma 3.6).  One does not need to deal with the set
up concerning the $f_i$ (they are all 0) and the estimates are similar to 
those of Lemma 6.2.  We omit a rigorous proof of Lemma 6.6 here and leave 
it to the reader.  We now run an argument similar to that which followed 
Lemma 6.2.  Dividing the dominated term above by $|\log \epsilon|$ in the 
conclusion of the statement of Lemma 6.4 and taking a lim sup as $\epsilon$ 
goes to $0$ yields $\beta_1^2 \cdot \bigtriangleup_{\varphi_1}(M_1)$.  Using 
the same asymptotic freeness results in the paragraphs preceding Theorem 
6.4 it follows that $\delta_0(a_1^{\prime},a_2^{\prime})$ is greater than or 
equal to this limiting process, i.e,  greater than or equal to $\beta_1^2 
\cdot \bigtriangleup_{\varphi_1}(M_1)$.  Since $M$ is hyperfinite ($M$ is 
a direct summand of matrix algebras), by Theorem 4.5 $\delta_0(a_1,\ldots,a_n) 
\geq \delta_0(a_1^{\prime}, a_2^{\prime}) \geq b_1^2 - \varepsilon - \sum_{i=1}
^{l} \frac {\alpha_i^2}{k_i^2}$.  $l$ and $\varepsilon$ being arbitrary 
$\delta_0(a_1,\ldots,a_n) \geq 1 - \sum_{i=1}^s \frac {\alpha_i^2}{k_i^2}.$    

	Secondly suppose $\alpha_0 = 1$.  By Section 4 
$\delta_0(a_1,\ldots,a_n) =1$ which yields the desired lower bound.  

	Having considered all the cases above we have for any $M$ and 
$a_1, \ldots, a_n$ as in Section 2:

\begin{theorem} $\delta_0(a_1,\ldots,a_n) \geq 1 - \sum_{i=1}^{s}
\frac{\alpha_{i}^{2}}{k_{i}^{2}}.$ \end{theorem}

By Theorem 3.9 we also have: \begin{corollary}If $M$ is hyperfinite,
then $\delta_0(a_1,\ldots,a_n)= 1 - \sum_{i=1}^{s}
\frac{\alpha_{i}^{2}}{k_{i}^{2}}$. \end{corollary}

	In light of Corollary 6.5 if $M$ is hyperfinite, then we define 
$\delta_0(M) = \delta_0(a_1,\ldots,a_n).$  As in the introduction we 
remark that every such hyperfinite $M$ has a finite set of self-adjoint 
generators.  

\section{Trivialities and a Final Remark}

	In concluding the discussion we make a few simple observations 
about the preceding results.  The first is a strengthening of Theorem 
4.5.  We start with a generalization of Lemma 3.8. 

\begin{corollary}If $N \subset M$ is a unital inclusion of hyperfinite
von Neumann algebras, then $\delta_0(N) \leq \delta_0(M)$.
\end{corollary} \begin{proof}Find self-adjoint generators
$a_1,\ldots,a_n$ for $M$ and $b_1,\ldots,b_p$ for $N$.  By Theorem 4.5
and Corollary 6.8,

\[ \delta_0(N) = \delta_0(b_1,\ldots,b_p) \leq 
\delta_0(a_1,\ldots,a_n,b_1,\ldots,b_p) = 
\delta_0(M). \]
\end{proof}

We now have:

\begin{corollary}(Hyperfinite Monotonicity).  If $N \subset M$ is a unital 
inclusion of von Neumann algebras and $N$ is hyperfinite, then $\delta_0(N) 
\leq \delta_0(a_1, \ldots,a_n).$
\end{corollary}
\begin{proof} We may assume $M = M_0 \oplus (\oplus_{i=1}^{s} M_{k_i}
(\mathbb C))$ and $\varphi = \alpha_0 \varphi_0 \oplus (\oplus_{i=1}^{s} 
\alpha_i tr_{k_i})$ where all quantities are as in Section 6.  
Define $A = \mathbb C I_0 \oplus (\oplus_{i=1}^{s} M_{k_i}(\mathbb C)) 
\subset M.$  It is easy to see that the von Neumann algebra $\mathcal R$ 
generated by $A \bigcup N$ is hyperfinite.  By decomposing $\mathcal R$ 
over its center and observing that the atomic projections of $Z(\mathcal R)$ 
contains those of $Z(M)$, it follows from Theorem 6.7 and Corollary 6.8 that $\delta_0
(\mathcal R) \leq \delta_0(a_1,\ldots,a_n).$  Hence by Corollary 7.1 $\delta_0(N) 
\leq \delta_0(\mathcal R) \leq \delta_0(a_1,\ldots,a_n).$
\end{proof}

	Our second observation is a weak lower semicontinuity property for 
$\delta_0$.  

\begin{lemma}If $\langle (b_{1}^{(k)},\ldots,b_{n}^{(k)})
\rangle_{k=1}^{\infty}$ is a sequence of n-tuples of self-adjoint
elements in $M$ such that for each $1 \leq i \leq n$ $b_{i}^{(k)}
\rightarrow b_i$ strongly and $\{b_1,\ldots,b_n\}$ generates a diffuse
von Neumann algebra, then

\[ \liminf_{k \rightarrow \infty} 
\delta_0(b_{1}^{(k)},\ldots,b_{n}^{(k)}) 
\geq 1. \]

\noindent In particular, if $1 = \delta_0(b_1,\ldots,b_n)$, then 
$\liminf_{k \rightarrow \infty}\delta_0(b_{1}^{(k)},\ldots,b_{n}^{(k)}) \geq 
\delta_0(b_1,\ldots,b_n)$.
 
\end{lemma}

\begin{proof} Suppose $\varepsilon >0$.  By the proof of Corollary 4.7 it 
follows that there exists a $b=b^* \in W^*(b_1,\ldots,b_n)$ such that 
$\delta_0(b) =1$.  There exists a sequence of noncommuting polynomials 
in $n$ variables $\langle q_m \rangle_{m=1}^{\infty}$ such that 
$q_m(b_1, \ldots, b_n)^* = q_m(b_1,\ldots,b_n) 
\rightarrow b$ strongly.  It follows from 
Proposition 6.14 of [9] and Corollary 6.7 of [10] that for $m$ 
sufficiently large $\delta_0(q_m(b_1,\ldots,b_n)) > 1- 
\frac{\varepsilon}{2}$.  Pick one such $m$ and call it $m_0$.  The same 
proposition of [9] and 
corollary of [10] provide a corresponding $N$ such that for all $k >N$  
$\delta_0(q_{m_0}(b_1^{(k)},\ldots,b_n^{(k)})) > 1 - \varepsilon$.  By 
Corollary 4.6 for all $k>N$ $\delta_0(b_{1}^{(k)},\ldots,b_{n}^{(k)}) \geq 
\delta_0(q_{m_0}(b_{1}^{(k)},\ldots,b_{n}^{(k)})) > 1- \varepsilon$.
\end{proof} 

	Finally, we comment on the work carried out by 
Ken Dykema concerning free products of hyperfinite von Neumann algebras 
with tracial, faithful states.  In [2] Dykema investigated the free 
product of two such algebras $A$ and $B$.  There it was shown that $A*B$ 
was isomorphic to $L(F_s) \oplus C$ where $L(F_s)$ is an interpolated 
free group factor and $C$ is a finite dimensional von Neumann algebra.  
Moreover, Dykema provided formulas for determining $C$ in terms of 
the matricial parts of $A$ and $B$ and calculating $s$ in terms of the 
the `free dimensions' of $A$, $B$, and $C$.  Given a hyperfinite $M$ 
as above, Dykema defined the free dimension of $M$, $\text {fdim} (M)$ to 
be 
\[ \alpha_0^2 + \sum_{i=1}^{s} \alpha_i^2 (1- k_{i}^{-2}) + 
2\alpha_0(1-\alpha_0) + \sum_{1 \leq i, j \leq s, i\neq j} \alpha_i 
\alpha_j.\]
Using the identity $1 = (\sum_{i=0}^{s} \alpha_i)^2$ one finds that the 
number above equals $\delta_0(M)$.  In other words, for a hyperfinite von Neumann 
algebra $M$ with a tracial, faithful state, the quantity $\delta_0(M)$ equals  
the quantity $\text{fdim} (M)$.

\section{Addendum}

	In this final section we prove the metric 
entropy estimates of Lemmas 3.5 and 5.2.  The proofs are essentially 
those of Szarek ([7]) with the addition of the explicit computations of 
Raymond ([5]).  

	Throughout $H$ will denote a closed Lie subgroup of $U_k$.  Define $X=U_k/H$, $|\cdot|_{\infty}$ to be 
the operator norm, $\mathcal H$ to be the Lie subalgebra of $H$ identified 
in $iM^{sa}_k(\mathbb C) = \mathcal G$, and $\mathcal X$ to be the 
orthogonal complement of $\mathcal H$ with respect to the real inner product on 
$\mathcal G$ generated by $Re$ $Tr$.  Denote $d_{\infty}$ 
and $d_2$ to be the metrics on $X$ induced by $| \cdot |_{\infty}$ and 
$| \cdot |_2$, respectively.  Lastly for a metric $d$ on a space 
$\Omega$ and $\epsilon >0$ define $N(\Omega,d,\epsilon)$ to be the minimum 
number of open $\epsilon$-balls required to cover $\Omega$ with respect to 
$d$ and $P(\Omega,d,\epsilon)$ to be the maximum number of points in 
an $\epsilon$-separated subset of $\Omega$ with respect to $d$.

	Szarek uses two essential quantities to obtain the 
metric entropy estimates in [7].  The first is $\kappa(M)$, 
the operator norm of the orthogonal projection onto $\mathcal X$ where the 
domain and range of the projection are equipped with 
the operator norm.  The second quantity Szarek employs are the weaving 
numbers of $X$.  We will use a slightly modified version of this.  The change is based 
on Szarek's preference to use the geodesic metric on $X$ and my 
inclination to use the extrinsic norm metric.  They are the same for our purposes.  
Given $\theta >0$ $H$ is $(\theta, | \cdot |_{\infty})$ 
-woven if for $u \in H, |u-I|_{\infty} < \theta \Rightarrow \exists h \in 
\mathcal H \hspace{.08in}\text {such that} \hspace{.08in} |h|_{\infty}< 
\frac{\pi}{16} \hspace{.08in}\text {and} \hspace{.08in} e^h = 
u$.  We define $\theta(X)$ to be the supremum over all $\theta$ satisfying 
the preceding condition. 

	We now state the main result of Szarek's ([7]), slightly altered
in our new notation.

\begin{theorem}Suppose $\beta \in (0, 1/2]$ and 
$\min \{ \theta (X), \kappa (X)^{-1}\} \geq \beta$.  Assume that one 
of the following conditions hold:

\begin{itemize}
\item $\dim H \leq (1-\beta)k^2$.
\item There exists a subspace $E \subset \mathbb C^k$ invariant under $H$ 
with $\dim E > \beta k$ satisfying $\beta k \leq dim E \leq 
(1 - \beta)k$.
\item There exists a subspace $E \subset \mathbb C^k$ invariant 
under $H$ with $p=\dim E > \beta k$ such that the decomposition $\mathbb 
C^k = E \oplus E^{\bot}$ induces an 
isomorphism $H \rightarrow U(p) \times H_o$ for some subgroup $H_o$ of 
$U_{k-p}$.
\end{itemize}
Then for any $\epsilon \in (0, \beta/4)$
\[ \left (\frac{c}{\epsilon} \right )^{\dim X} \leq 
N(X,d_{\infty},\epsilon) \leq \left (\frac{C}{\epsilon} \right )^{\dim X} 
\]
where $c, C>0$ are constants depending only on $\beta$. 
\end{theorem} 

The utility of Szarek's result lies in the fact that the quantities $c$
and $C$ depend only on $\beta$.  We now provide the proof of Lemma 3.5.  

\noindent\emph{Proof of Lemma 3.6.}\hspace{.1in}Suppose $H$ is tractable.
Consider the conditional expectation $e : M_k(\mathbb C)  \rightarrow
H^{\prime \prime}$.  $I - e$ restricted to $\mathcal G$ is the orthogonal
projection onto $\mathcal X$ and since $\|e\| \leq 1$, it follows that
$\kappa (X)^{-1} \geq \frac{1}{2}$.  The spectral theorem shows that
$\theta (X) > |e^{i \frac{\pi}{16}} - 1|$.  Hence, $\min \{\theta (X),
\kappa (X)^{-1} \} > \frac {1}{20}$.  I claim that $H$ satisfies one of
the three conditions as stated in the theorem for $\beta = \frac {1}{20}$.  
Without loss of generality assume $H$ takes the form appearing
in the definition of a tractable Lie subgroup of $U_k.$
Suppose there exist some $1 \leq j_1, \ldots, j_q \leq m$ such that $
\frac {k}{20} \leq \sum_{i=1}^{q} k_{j_i}l_{j_i} \leq \frac{19k}{20} $.  
Then $H$ satisfies the second condition of Theorem 8.1.  Otherwise there
must exist some $1 \leq i \leq m$ for which $k_il_i > \frac {19k}{20}$.  
If $k_i =1$, then $H$ satisfies the third condition of the theorem.  
Otherwise $k_i >1$ and this forces there to be a reducing subspace $E$ for
$H$ with $ \frac{19k}{60} \leq \dim E \leq \frac {k}{2}$ whence $H$
fulfills the second condition of the theorem.  Theorem 8.1 now yields the
desired result.  \hspace{3.7in}$\square$

	Having dealt with Lemma 3.5 let's turn to the finite dimensional
situation in Lemma 5.2.  More generally first consider the viability of
the lower bounds of Theorem 8.1 when $X$ is obtained from tractable $H$
and where instead of using $d_{\infty}$ we use $d_2.$ Some results of [7]
works for unitarily invariant norms and metrics but with the $| \cdot
|_2$-norm problems arise.  The quantity $\theta(X)$, properly interpreted
does not stay uniformly away from $0$ even when we consider the
homogeneous spaces in Section 5 associated to a finite dimensional $M$.  
Presumably $X$ would be $(\theta, | \cdot|_2)$-woven if for $u \in H, |u -
I|_2 < \theta \Rightarrow \exists h \in \mathcal H$ such that
$|h|_{\infty} < \frac {\pi}{16}$ and $e^h=u$.  Unfortunately, the
homogeneous spaces which we will restrict our attention to (which is much
smaller than the class of homogeneous spaces obtained from tractable $H$)
will fail to have $\theta$ values uniformly bounded away from $0$.  
Nevertheless, we still have the key result [7], Lemma 10, where the use of
$\theta(X)$ was crucial:

\begin{lemma}There exist $\lambda, r >0$ such that for any $k \in
\mathbb N$ and tractable $H$ of $U_k$ if $x,y \in \mathcal
X$, and $|x|_{\infty},|y|_{\infty} < r$, then

\[ d_2(q(e^x),q(e^y)) \geq \lambda |x-y|_2 \]

\noindent where $q: U_k \rightarrow X$ is the quotient map. \end{lemma}
\begin{proof}For $r$ (as yet to be specified) and any such $x$ and $y$ as
above, there exists by definition an $h \in \mathcal H$ with $|h|_{\infty}
\leq \pi$ satisfying $d_2(q(e^x),q(e^y)) = \inf_{v \in H} |e^{-y}e^x-v|_2
= |e^{-y}e^x - e^h|_2$. Set $u = e^{-y}e^x$.  By the spectral theorem
write $ h = i \sum_{j=1}^{d} \beta_j f_j$ where the $f_j$ are mutually
orthogonal projections and the $\beta_j$ are real numbers.  We can arrange
it so that for each $j$, $i f_j \in \mathcal H$, i.e., $h$ takes the block
form of $\mathcal H$.  Define $\gamma_j$ to be $4r$ if $\beta_j > 4r$,
$-4r$ if $\beta_j < -4r$, and $\beta_j$ if $|\beta_j| \leq 4r$.  Set $z =
i \sum_{j=1}^{d} \gamma_j f_j \in \mathcal H.$ $|z|_{\infty} \leq 4r.$
Define $\Lambda_1 = \{ j \in \mathbb N : 1 \leq j \leq d, |\beta_j| \leq
4r\}$ and $\Lambda_2 = \{1,\ldots,d\} - \Lambda_1$.  Observe that
$|x|_{\infty}, |y|_{\infty} < r \Rightarrow |u - I|_{\infty} < 2r.$

\begin{eqnarray*} |u- e^z|_2^2 \leq \sum_{j=1}^d |uf_j -
e^{i\gamma_j}f_j|_2^2 &=& \sum_{j\in \Lambda_1} |uf_j -
e^{i\beta_j}f_j|_2^2 + \sum_{j\in \Lambda_2} |uf_j- e^{i\gamma_j}f_j|_2^2
\\ & \leq & |u - e^h|_2^2 + \sum_{j \in \Lambda_2} (6r)^2 |f_j|_2^2.  
\end{eqnarray*}

\noindent Now for $r$ sufficiently small $3r \leq |1 - e^{i4r}|$ ($r$ 
dependent only upon the exponential map).  

\begin{eqnarray*} \sum_{j \in \Lambda_2} (6r)^2|f_j|_2^2 \leq \sum_{ j \in
\Lambda_2} 36 ( | 1 - e^{i4r} | - 2r)^2 |f_j|_2^2 & \leq & \sum_{j \in
\Lambda_2} 36 (|1 - e^{i\beta_j}| - 2r)^2 |f_j|_2^2 \\ & \leq & \sum_{j
\in \Lambda_2} 36 ( |f_j - e^{i\beta_j}f_j|_2 - |f_j -uf_j|_2 )^2 \\ &
\leq & \sum_{j \in \Lambda_2} 36|uf_j - e^{i\beta_j}f_j|_2^2 \\ & \leq &
36|u-e^h|_2^2.  \end{eqnarray*}

\noindent It follows that $|u - e^z|_2 \leq 7|u-e^h|_2.$  

	A repetition of the proof of Lemma 10 in [7] minus the parts
referring to $\theta(X)$ shows that there exist $\lambda, r >0$ (which we
can make as small as we want and in particular have $r$ satisfy $3r \leq
|1-e^{i4r}|$) independent of the tractable $H$ such that for any $x,y\in
\mathcal X$ with $|x|_{\infty}, |y|_{\infty} < r$, $|e^{-y}e^x - e^z|_2
\geq \lambda |x-y|_2$. By what preceded for any $x, y \in \mathcal X$ with
$|x|_{\infty}, |y|_{\infty} < r$

\[d_2(q(e^x),q(e^y)) = |e^{-y}e^x - e^h|_2 \geq \frac{|e^{-y}e^x - 
e^z|_2}{7} \geq \frac {\lambda |x-y|_2}{7}.\]
\end{proof}

The result above does not quite provide the desired lower bounds for the
homogeneous space associated to $H.$ Observe that if all the normalized
Hilbert-Schmidt quantities are replaced with operator norm quantities,
then $P(X,d_{\infty},\epsilon)$ is bounded below by the
$\frac{\epsilon}{\lambda}$ packing number of the $r$-ball of the space
$\mathcal X$ (endowed with $|\cdot|_{\infty}$).  The appropriate lower
bounds for packing numbers of balls in finite dimensional spaces can be
obtained through a standard volume comparison argument (see [1] for an
example of how this technique yields the packing number bounds).  Indeed,
this is how the lower bound is achieved in Theorem 8.1.  The result above
is not quite the same.  It shows that $P(X,d_2,\epsilon)$ dominates the
$\frac{\epsilon}{\lambda}$ packing number with respect to the $| \cdot
|_2$-metric of the ball of $| \cdot |_{\infty}$-radius r in $\mathcal X.$
The issue is that we have a lower bound involving \emph{two} different
metrics.  We want to obtain the appropriate lower bounds by using the
volume comparison argument but our task is slightly complicated by this.  
We must now examine the ratio of the volumes of balls of
radius 1 with respect to $| \cdot |_{\infty}$ and $|\cdot|_2$ in the space
$\mathcal X$ associated to $H.$

	Despite the difficulties mentioned Lemma 5.2 demands 
lower bounds on the packing numbers of a specific class of homogeneous 
spaces, in fact, much smaller than the class of all homogeneous spaces 
obtained from tractable subgroups.  Hence, the task at hand is not so 
daunting.  With Lemma 8.2 in hand we now begin the main part of the 
proof of Lemma 5.2.  As discussed at the end of the preceding paragraph, 
our main objective is to examine the ratio of 
the volumes of balls in the orthogonal complements of certain Lie subalgebras.  

\noindent\emph{Proof of Lemma 5.2.} \hspace{.03in} Maintain all the
assumptions made on $M$ and $\varphi$ in Section 5.  We assume that $M
\neq \mathbb C I$ since Lemma 8.2 clearly holds in this situation.  
There exist constants $1 > \delta, c_1, c_2 >0$ such that if $\delta >
\varepsilon >0$ and $r_1,\ldots,r_p \in \mathbb R$ satisfy $|r_j - \frac
{\alpha_j}{n_j}| < \varepsilon$ for all $j,$ then $c_1 < \sum_{j=1}^p
r_j^2 < c_2$ ($M \neq \mathbb C I \Rightarrow \sum_{j=1}^p \left (\frac
{\alpha_j}{n_j} \right)^2 < 1$). Now suppose $ \frac {1}{2} \cdot \min
\{ \delta, 1-c_1, \alpha_1, \ldots, \alpha_p \} > \varepsilon >0.$ It is
a trivial consequence of Lemma 3.6 that there exists a sequence $\langle
\sigma_k \rangle_{k=1}^{\infty}$ such that for each $k$ $\sigma_k:M
\rightarrow M_k(\mathbb C)$ is a $*$-homomorphism and for $k$
sufficiently large: \begin{itemize} \item $\|tr_k \circ \sigma_k -
\varphi \| < \varepsilon$. \item The set of unitaries $H_k$ of
$\sigma_k(M)^{\prime}$ is a tractable Lie subgroup of $U_k$ and setting
$X_k = U_k/H_k$ we have that $k^2(\bigtriangleup_{\varphi}(M) -
\varepsilon) \leq \dim(X_k).$ \end{itemize} We must demonstrate the
third item in Lemma 5.2 (the lower bound packing estimate) and make sure
that the constant $\kappa$ obtained is independent of $\varepsilon$ and
$\epsilon.$ It can be arranged so that there exists a $k_0 \in \mathbb
N$ such that for all $k \geq k_0$ the representation $\sigma_k$ takes
the simple form described in the proof of Lemma 3.6.  Recall from the
proof of Lemma 3.6 that for each $k \geq k_0$ we have the
$l_1(k),\ldots,l_{p+1}(k),$ associated to $\sigma_k.$ Moreover for such
$k$ the construction of the $\sigma_k$ in Lemma 3.6 and the bound placed
on $\varepsilon$ shows that $ \frac {3}{2} \alpha_j> \frac
{l_j(k)n_j}{k} > \frac {\alpha_j}{2}$ for each $j,$ $ c_1 k^2 <
\sum_{j=1}^{p+1} l_j(k)^2 < c_2 k^2,$ and $l_{p+1}(k) < n_1 \cdots n_p$
for $k > k_0.$

	For each $k$ denote $\mathcal H_k$ and $\mathcal X_k$ to be the 
spaces $\mathcal H$ and $\mathcal X$ associated to $H=H_k$.  Again, we've 
translated the packing number problem (the third condition of Lemma 5.2) 
into the problem of comparing the volumes of the balls of $\mathcal X_k$ 
with respect to the norms $| \cdot |_{\infty}$ and $| \cdot |_2$ and finding 
an appropriate relationship between the two values for sufficiently large $k.$

	For any $r>0$ denote by $\mathcal G_k^r, \mathcal H_k^r$, and
$\mathcal X_k^r$ the balls centered at the origin of operator norm less
than or equal to $r$ in $iM_{k}^{sa}(\mathbb C), \mathcal H_k$, and
$\mathcal X_k$, respectively.  Consider the conditional expectation $e$
for $H_k^{\prime \prime}$.  Define $ \Phi :  iM_{k}^{sa}(\mathbb C)
\rightarrow \mathcal H_k \oplus \mathcal X_k$ by $\Phi(x) = (e(x),
(I-e)x)$.  Since $e$ is a contraction when its domain and range are
endowed with the operator norm, it follows that $\Phi(\mathcal G_k^1)  
\subset \mathcal H_k^1 \times \mathcal X_k^2$.  $\Phi$ is an isometry when
its domain and range are endowed with the Hilbert-Schmidt norm (normalized
or not).  Thus, $\text{vol} (\mathcal G_k^1) = \text{vol} (\Phi(\mathcal
G_k^1)) \leq \text{vol} (\mathcal H_k^1) \cdot \text{vol} (\mathcal
X_k^2)$.  Notice that here we calculate the volumes of $\mathcal H_k^1$
and $\mathcal X_k^2$ in their ambient Hilbert spaces $\mathcal H_k$ and
$\mathcal X_k$ endowed with the real inner product $Re$ $Tr.$

We claim that if for each $d$ $\Lambda_d = \text{vol}(\mathcal G^1_d)$ and
$\Theta_d$ denotes the volume of the ball of radius $\sqrt{d}$ in $\mathbb
R^{d^2},$ then there exists some constant $1> \zeta_1 >0$ such that
$(\zeta_1)^{d^2} \leq \frac {\Lambda_d}{\Theta_d}$ for all $d.$ By [5]
there exists a constant $c>0$ such that $\frac {a_d}{b_d} \sim (c)^{2d^2}$
as $d \rightarrow \infty$ where $a_d$ is the volume (with respect to the
real inner product $Re$ $Tr$) of the operator norm unit ball of
$M_k(\mathbb C)$ and $b_d$ is the volume of the ball of radius $\sqrt{2d}$
in $2d^2$-dimensional real Euclidean space.  We now use the same trick in
the preceding paragraph.  Decompose $M_k(\mathbb C)$ as the orthogonal
direct sum $M_{k}^{sa}(\mathbb C) \oplus iM_{k}^{sa}(\mathbb C)$.  It
follows that the operator norm unit ball of $M_k(\mathbb C)$ is contained
in the direct sum of the operator norm unit ball of $M_{k}^{sa}(\mathbb
C)$ and the operator norm unit ball of $\mathcal G_d^1$.  The volume of
this latter set is $(\Lambda_d)^2.$ So $ \frac {a_d}{b_d} < \frac
{(\Lambda_d)^2}{b_d} \sim \frac {(\Lambda_d)^2}{(\Theta_d)^2}$ and this
yields the desired result.

	Now $\text{vol} (\mathcal H_k^1) <
\Theta_{l_1(k)}\cdots\Theta_{l_{p+1}(k)} \cdot \sqrt {n_1 \cdots n_p}$ and
because $\dim \mathcal X_k > (1- c_2 ) \cdot k^2$ for
sufficiently large $k,$ it follows that there exists a $\zeta_2 >0$
(dependent only on $n_1,\ldots,n_p, \zeta_1$ and $c_1$) such that for
sufficiently large $k$

\[ \text{vol} (\mathcal X_k^2) \geq \frac {\text{vol} (\mathcal G_k^1)}{
\text{vol} (\mathcal H_k^1)} \geq \frac {\Theta_k}{\Theta_{l_1(k)} \cdots
\Theta_{l_{p+1}(k)}} \cdot \zeta_2^{\dim \mathcal X_k}. \]

\noindent $\frac {l_j(k)}{k} > \frac {\alpha_j}{2n_j},$ $l_{p+1}(k) < n_1
\cdots n_p$ for $k > k_0,$ and $(1- c_2) \cdot k^2 < \dim \mathcal X_k <
(1 -c_1) \cdot k^2$ for $k > k_0.$ These three facts and Stirling's
formula shows that there exists a constant $\zeta_3 >0$ (dependent only on
the $\alpha_i, n_i, c_1,$ and $c_2$) such that for sufficiently large $k$
the dominated term above is greater than or equal to $ (\zeta_3)^{\dim
\mathcal X_k} \cdot C_{\dim \mathcal X_k}$ where $C_{\dim \mathcal X_k}$
is the volume of the ball of radius $\sqrt{k}$ in $\mathbb R^{\dim
\mathcal X_k}.$ Notice that this quantity is the volume of the ball of
$\mathcal X_k$ of $| \cdot|_2$ radius 1.  In other words there exists a
$\zeta >0$ (again dependent only on the $\alpha_i, n_i, c_1,$ and $c_2$)
such that for sufficiently large $k$

\[ \frac{ \text{vol} (\mathcal X_k^1)}{C_{\dim \mathcal X_k}} > (\zeta 
)^{\dim 
\mathcal X_k}. \]  

The standard volume comparison method (for an example of how this method
is used see [1]) shows that for such $k$ and any $\epsilon>0$ $P(\mathcal
X_k^r, |\cdot|_2,\epsilon) > \left (\frac {r \zeta}{ \epsilon}
\right)^{\dim \mathcal X_k}$.  Using Lemma 8.2 it follows that for $k$
sufficently large, if $\epsilon >0$, then

\[ P(X_k, d_2, \epsilon) > \left (\frac {\lambda r \zeta}{ \epsilon} \right 
)^{\dim X_k}.\]
       
\noindent Set $\kappa =  \lambda r \zeta.$ $\kappa$ depends
only on $\alpha_i, n_i, c_1, c_2,$ and the upper bound placed on
$\varepsilon.$ The upper bound can be relaxed for the purposes of
Lemma 5.2.  We have the third and final condition of Lemma 5.2.
\hspace{6.25in}$\square$ \vspace{.1in}

\noindent{\it Acknowledgements.} I thank Dan Voiculescu, my advisor, both
for suggesting that I work on computing the free entropy dimension of
$M_2(\mathbb C)$ and for providing helpful comments throughout.  Also I
thank Stanislaw Szarek for the exchanges and suggestions which formed the
core of the proof of Lemma 5.2 in the Addendum.  I am grateful to Charles
Holton for the simplification of the proof of Corollary 7.1 and the useful
conversations on Minkowski dimension.  Finally, this work would not have
been possible without the generous contributions of the NSF Graduate
Fellowship program.


\begin{thebibliography}{[ASMR]} 

\bibitem{1}Carl, B. and Stephani, I. {\it Entropy , Compactness, and the 
Approximation of Operators}, Cambridge University Press, Cambridge, 1990.

\bibitem{2} Dykema, K. J., {\it Free Products of hyperfinite von Neumann 
algebras and free dimension} Duke Math. J. 69 (1993), 97-119.

\bibitem {3} Ge, Liming {\it Applications of free entropy to finite von 
Neumann algebras, II}, Annals of Mathematics, 147 (1998), 143-157.

\bibitem{4} Ge, Liming and Shen, Junhao {\it On Free Entropy Dimension of
Finite Von Neumann Algebras}, Geometric and Functional Analysis, Vol. 12, 
(2002), 546-566.

\bibitem{5} Raymond, Jean Saint, {\it Le volume des id\'{e}aux 
d'op\'{e}rateurs classiques}, Studia Mathematica, v.LXXX (1984), 63-75. 

\bibitem {6} Stefan, M., {\it The indecomposability of free group factors 
over nonprime subfactors and abelian subalgebras}, preprint.
  
\bibitem{7} Szarek, S. {\it Metric Entropy of homogeneous spaces}, Quantum 
Probability, (Gdensk, 1997), Banach Center Publications v.43, 
Polish Academy of Science, Warsaw 1998, 395-410.

\bibitem{8} Voiculescu, D., Dykema, K.J., Nica, A. {\it Free Random 
Variables}, CRM Monograph Series, v. 1, American Mathematical Society, 
1992.

\bibitem{9} Voiculescu, D. {\it The Analogues of entropy and of Fisher's 
information measure in free probability theory, II}, Inventiones 
mathematicae 118, (1994), 411-440.

\bibitem{10} Voiculescu, D. {\it The Analogues of Entropy and of Fisher's 
Information Measure in Free Probability Theory III: The Absence of Cartan 
Subalgebras}, Geometric and Functional Analysis, Vol. 6, No.1 (1996) 
(172-199).

\bibitem{11} Voiculescu, D. {\it A Strengthened Asymptotic Freeness Result
for Random Matrices with Applications to Free Entropy}, IMRN, 1 (1998),
41-64.

\end{thebibliography}
\end{document}